\RequirePackage{fix-cm}

\RequirePackage{amsmath}
\documentclass[smallcondensed]{svjour3}     
\smartqed  
\usepackage{graphicx}
\usepackage{subfigure}
\usepackage{algpseudocode}
\usepackage{algorithmicx,algorithm}

\usepackage{graphicx}
\usepackage{amsfonts}
\usepackage{amssymb}
\usepackage{authblk}

\usepackage{color}

\usepackage{natbib}
\bibpunct{(}{)}{;}{}{}{;}

\begin{document}

\title{Asynchronous Stochastic Composition Optimization with Variance Reduction
}


\author{Shuheng Shen         \and
        Linli Xu             \and
        Jingchang Liu        \and
        Junliang Guo         \and
        Qing Ling       
}


\institute{Shuheng Shen \at
           University of Science and Technology of China,~Hefei,~China \\
              \email{vaip@mail.ustc.edu.cn}           
           \and
           Linli Xu \at
           University of Science and Technology of China,~Hefei,~China \\
              \email{linlixu@ustc.edu.cn}
           \and
           Jingchang Liu \at
           University of Science and Technology of China,~Hefei,~China \\
              \email{xdjcl@mail.ustc.edu.cn}
           \and
           Junliang Guo \at
           University of Science and Technology of China,~Hefei,~China \\
              \email{leoguojl@gmail.com}
           \and
           Qing Ling \at
           Sun Yat-Sen University,~Guangzhou,~China \\
              \email{lingqing556@mail.sysu.edu.cn}
}

\date{Received: date / Accepted: date}

\maketitle

\begin{abstract}
Composition optimization has drawn a lot of attention in a wide variety of machine learning domains from risk management to reinforcement learning. Existing methods solving the composition optimization problem often work in a sequential and single-machine manner, which limits their applications in large-scale problems. To address this issue, this paper proposes two asynchronous parallel variance reduced stochastic compositional gradient (AsyVRSC) algorithms that are suitable to handle large-scale data sets. The two algorithms are AsyVRSC-Shared for the shared-memory architecture and AsyVRSC-Distributed for the master-worker architecture. The embedded variance reduction techniques enable the algorithms to achieve linear convergence rates. Furthermore, AsyVRSC-Shared and AsyVRSC-Distributed enjoy provable linear speedup, when the time delays are bounded by the data dimensionality or the sparsity ratio of the partial gradients, respectively. Extensive experiments are conducted to verify the effectiveness of the proposed algorithms. 


\end{abstract}

\section{Introduction}
\label{intro}
Consider the problem of composition optimization~\citep{wang2017stochastic} which minimizes a loss function with a compositional expected form:
{\small
\begin{equation} \label{expectation_formula}
\min_{x \in \mathbb{R}^{d_1}} f(x):=\mathbb{E}_{i} F_{i}(\mathbb{E}_{j} G_{j}(x)),
\end{equation}
}where $G_{j}(x): x \in \mathbb{R}^{d_1} \mapsto y \in \mathbb{R}^{d_2}$ are inner component functions and $F_{i}(y): y \in \mathbb{R}^{d_2} \mapsto z \in \mathbb{R}$ are outer component functions, both of which are continuously differentiable. Many emerging applications can be formulated as problem~(\ref{expectation_formula}), such as reinforcement learning~\citep{dai2016learning}, risk-reverse learning~\citep{wang2016accelerating}, multi-stage stochastic programming~\citep{shapiro2009lectures}, adaptive simulation~\citep{hu2014model}, etc.

In practice, the expectation in (\ref{expectation_formula}) can be replaced by a finite-sum form when the number of samples is finite. This paper focuses on the finite-sum composition optimization problem:
{\small
\begin{equation} \label{basic_formula}
\min_{x \in \mathbb{R}^{d_{1}}} f(x):=\frac{1}{n_1} \sum_{i=1}^{n_1} F_{i}( \frac{1}{n_2} \sum_{j=1}^{n_2} G_{j}(x)),
\end{equation}
}where $n_1$ is the number of outer samples and $n_2$ is the number of inner samples. For the ease of presentation, we use $G(x):=\frac{1}{n_2} \sum_{j=1}^{n_2} G_{j}(x)$ and $F(G(x)):=\frac{1}{n_1} \sum_{i=1}^{n_1} F_{i}(G(x))$ to denote the inner function and the outer function, respectively. Then the full gradient of $f(x)$ can be represented as $\nabla f(x) = (\nabla G(x))^{T} \nabla F(G(x))$, where $\nabla G(x) \in \mathbb{R}^{d_2 \times d_1}$ is the Jacobian of $G(x)$.

Composition optimization is substantially more challenging than a general optimization problem in the form of $\min_{x \in \mathbb{R}^{d_1}} \frac{1}{n_1} \sum_{i=1}^{n_1} F_{i}(x)$. Standard stochastic gradient descent (SGD) algorithm is not well suited for minimizing composition problems since it needs to calculate $\nabla G_{j}(x) \nabla F_{i}(G(x))$ at each iteration, which is time-consuming for the computation of a full inner function value $G(x)$. To address this issue, stochastic compositional gradient descent (SCGD) and its accelerated version are proposed  in~\citep{wang2017stochastic}; both of them have a constant query complexity per iteration. However, the variance introduced by random sampling in SCGD results in a sublinear convergence rate, even for strongly convex loss functions. This fact motivates the combination of variance reduction techniques, which have been successfully applied in SGD~\citep{johnson2013accelerating,defazio2014saga,schmidt2017minimizing}, with SCGD~\citep{lian2017finite}. Specifically, along the line of stochastic variance reduced gradient (SVRG)~\citep{johnson2013accelerating}, the authors of \citep{lian2017finite} propose two variance reduced stochastic composition algorithms, Composition-SVRG-1 and Composition-SVRG-2, both of which have provable linear convergence rates for strongly convex loss functions and the later performs much better if the loss function has a large condition number. Composition-SVRG-2~(denoted as VRSC in this paper) proposed in~\citep{lian2017finite} has the state-of-art performance for composition optimization in the single-machine setting.

In the meantime, with the growth of sample size and model complexity, it becomes challenging to train machine learning models on large-scale datasets using sequential~(single-machine) algorithms. Therefore, asynchronous parallel algorithms that distribute computation to multiple workers are popular solutions to the scalability issue. Successful applications include asynchronous parallel implementations of SGD~\citep{niu2011hogwild,lian2015asynchronous}, SVRG~\citep{reddi2015variance,meng2017asynchronous,huo2017asynchronous}, stochastic coordinate descent (SCD)~\citep{liu2015asynchronous,richtarik2016distributed}, SAGA~\citep{leblond2017asaga,pedregosa2017breaking}, etc. In these algorithms, the workers calculate sample gradients and update the optimization variable in a parallel and asynchronous manner. However, asynchronous parallelization of composition optimization remains an open problem, since the compositional structure brings essential difficulties to parallel computation.

This paper proposes two asynchronous parallel variance reduced stochastic composition (AsyVRSC) algorithms, which fit for large-scale applications, have linear convergence guarantee for strongly convex loss functions, and enjoy linear speedup with respect to the number of workers. To be specific, the algorithms are developed for two major distributed computation architectures, shared-memory for multi-core or multi-GPU systems~\citep{niu2011hogwild} and master-worker for multi-machine clusters~\citep{agarwal2011distributed}. We prove that the proposed algorithms have linear speedup when time delays are bounded by data dimensionality in the shared-memory architecture, or bounded by sparse ratio of partial gradients in the master-worker architecture, demonstrating their potential to solve large-scale problems in asynchronous parallel environments.

The contributions of our work are listed as follows:
\begin{itemize}
    \setlength{\itemsep}{0pt}
    \setlength{\parsep}{0pt}
    \setlength{\parskip}{0pt}
\item We propose two asynchronous parallel variance reduced stochastic compositional gradient algorithms, AsyVRSC-Shared and AsyVRSC-Distributed for the shared-memory and master-worker architectures, respectively. 
\item We prove that both AsyVRSC-Shared and AsyVRSC-Distributed can achieve linear speedup with respect to the number of workers under certain conditions. 
\item Experiments on two tasks including portfolio management and reinforcement learning verify the effectiveness of the proposed algorithms.
\end{itemize}

\noindent \textbf{Notations.} We use $\|x\|$ to denote the $L_2$-norm of $x$, and $\langle x,y \rangle$ to denote the inner product of $x$ and $y$. The set $\{1,2,\cdots,n\}$ is represented by $[n]$. $\nabla_{k} f_{i}(x)$ indicates the $k$-th coordinate of the vector $\nabla f_{i}(x)$ and $\nabla_{k,l} G_{j}(x)$ corresponds to the $(k,l)$-th entry of the matrix $\nabla G_{j}(x)$. We denote by $\mathbb{E}$ a full expectation with respect to all the randomness.

\section{Preliminary: VRSC}
\label{pre}
We first briefly review the variance reduced stochastic compositional gradient (VRSC) method that solves the finite-sum composition optimization problem (\ref{basic_formula})~\citep{lian2017finite}.

Similar to SVRG that is a variance reduced modification to SGD, VRSC~\citep{lian2017finite} has two loops. In the \emph{s-th} outer loop, given an initial point $\widetilde{x}^s$, one keeps a snapshot of $G(\widetilde{x}^s)$,$\nabla G(\widetilde{x}^s)$ and $\nabla f(\widetilde{x}^s)$:
{\small
\begin{equation}\label{out_G}
G(\widetilde{x}^s) = \frac{1}{n_2} \sum_{j=1}^{n_2} G_{j}(\widetilde{x}^s),
\end{equation}
\begin{equation}\label{out_delta_G}
\nabla G(\widetilde{x}^s) = \frac{1}{n_2} \sum_{j=1}^{n_2} \nabla G_j(\widetilde{x}^{s}),
\end{equation}
\begin{equation}\label{full_grad}
\nabla f(\widetilde{x}^s) = (\nabla G(\widetilde{x}^s))^T \nabla F(G(\widetilde{x}^s)),
\end{equation}
}where $G(\widetilde{x}^s)$ denotes the value of the inner function, $\nabla G(\widetilde{x}^s)$ denotes the gradient of the inner function and $\nabla f(\widetilde{x}^s)$ is the full gradient. At the \emph{t-th} iteration of the inner loop, to estimate the gradient $\nabla f(x_t^s)$ at the current parameter $x_{t}^{s}$, $G(x_{t}^{s})$ and $\nabla G(x_{t}^{s})$ are estimated first by uniformly sampling two mini-batches $A_{t}$ and $B_{t}$ from $[n_{2}]$ with size $a$ and $b$ respectively:
{\small
\begin{equation} \label{in_G}
\widehat{G}_{t}^{s} = G(\widetilde{x}^s) - \frac{1}{a} \sum_{j=1}^{a}( G_{A_{t}[j]}(\widetilde{x}^{s}) - G_{A_{t}[j]}(x_{t}^{s}) ),
\end{equation}
\begin{equation} \label{in_delta_G}
\nabla \widehat{G}_{t}^{s} = \nabla G( \widetilde{x}^{s} ) - \frac{1}{b} \sum_{j=1}^{b}( \nabla G_{B_{t}[j]}(\widetilde{x}^{s}) - \nabla G_{B_{t}[j]}(x_{t}^{s}) ),
\end{equation} 
}where $A_t[j]$ and $B_t[j]$ stand for the \emph{j-th} elements of $A_t$ and $B_t$, respectively. Based on the estimation of $G(x_{t}^{s})$ and $\nabla G(x_{t}^{s})$, $\nabla f(x_{t}^{s})$ can be estimated by:
{\small 
\begin{equation} \label{in_delta_f}
\nabla \widehat{f}(x_{t}^{s}) 
= (\nabla \widehat{G}_{t}^{s})^{T} \nabla F_{i_t}(\widehat{G}_{t}^{s}) - \nabla f_{i_t}(\widetilde{x}^{s}) + \nabla f( \widetilde{x}^{s} ),
\end{equation} 
}where $\nabla f_{i_t}(\widetilde{x}^{s}) = (\nabla G(\widetilde{x}^{s}))^T \nabla F_{i_t}(G(\widetilde{x}^{s}))$ and $i_{t}$ is uniformly sampled from $[n_1]$. This way, VRSC reduces the variance of SCGD, and improves the convergence rate from sublinear to linear under a constant learning rate for strongly convex problems~\citep{lian2017finite}.

VRSC has achieved great success for accelerating the minimization of composition optimization problems. However, it is still time-consuming when the data scale is large. To further accelerate the optimization, in this paper, we propose two asynchronous parallel algorithms for stochastic composition optimization with variance reduction. In the next sections, we introduce the two algorithms, followed by their theoretical analysis.

\section{Our Algorithms}
\label{alg}
In this section, we propose two asynchronous parallel variance reduced stochastic composition optimization algorithms, AsyVRSC-Shared that fits for the shared-memory architecture and AsyVRSC-Distributed for the master-worker architecture.

\begin{algorithm}[!htb]
\caption{AsyVRSC-Shared}
\label{alg:shared}
\hspace*{0.02in}{\bf Input:}
Inner iteration number $K$, outer iteration number $S$, mini-batch sizes $a$, $b$, learning rate $\eta$, initial point $\widetilde{x}^{1} \in \mathbb{R}^{d_1}$. 
\begin{algorithmic}[1]
\For{$s=1,2,...,S$} 
    \State \textbf{Phase 1:}
    \State Synchronously compute $G(\widetilde{x}^s)$, $\nabla G(\widetilde{x}^s)$ and $\nabla f(\widetilde{x}^s)$ using (\ref{out_G}), (\ref{out_delta_G}) and (\ref{full_grad}) 
    \State $x_{0}^{s} = \widetilde{x}^s$
    \State \textbf{Phase 2:}
    \For{$t=0,1,2,...,K-1$, \emph{asynchronously}}
        \State Read $x_{t-\tau_t^s}^{s}$ from the shared memory
        \State Uniformly sample $A_t$ and $B_t$ from $[n_2]$ with replacement, where $|A_t|=a, |B_t|=b$
        \State Compute $\nabla \widehat{f}(x_{t-\tau_t^s}^{s})$ using (\ref{in_G}), (\ref{in_delta_G}) and (\ref{in_delta_f})
        \State Uniformly sample $k_t$ from $[d_1]$
        \State Update $(x_{t+1}^{s})_{k_t} = (x_{t}^{s})_{k_t} - \eta (\nabla \widehat{f}(x_{t-\tau_t^s}^{s}))_{k_t}$
    \EndFor
    \State $\widetilde{x}^{s+1} = x_{r}^{s}$ for randomly chosen $r \in \{0,\cdots, K-1\} $
\EndFor
\State \Return $\widetilde{x}^{S+1}$
\end{algorithmic}
\end{algorithm}

\subsection{AsyVRSC-Shared}
In a shared-memory architecture, suppose there are $W$ local workers, each of which has full access to the whole training data and the parameters. Each local worker independently reads the parameter from the shared memory, computes a stochastic gradient and updates the parameter~\citep{niu2011hogwild}.

Note that VRSC does not have a vanilla asynchronous parallel implementation since it has two phases in each outer loop, which
makes it impossible to be completely asynchronous. We implement AsyVRSC-Shared summarized in Algorithm~\ref{alg:shared} to synchronously calculate the full gradient, and
to asynchronously update the parameter. 

\textbf{Phase $1$:} As line $3$ in Algorithm~\ref{alg:shared}, at the beginning of the \emph{s-th} outer loop, we keep a snapshot of $G(\widetilde{x}^s)$, $\nabla G(\widetilde{x}^s)$ and $\nabla f(\widetilde{x}^s)$ computed by all workers synchronously. 

\textbf{Phase $2$:} In the inner loops, all workers calculate the gradients and update the parameter in shared memory independently in an asynchronous way, corresponding to Algorithm~\ref{alg:shared}, line $6$-$12$.

AsyVRSC-Shared is a lock-free implementation, which means that the parameter in the shared memory may be updated while a worker is reading it. Therefore, the parameter one worker reads from the shared memory may be not a real state of $x$ at any time point. To avoid this inconsistency as much as possible, we calculate and update a single component (one coordinate) of the parameter since updating a single component of the parameter can be viewed as an atomic operation. This is along the line of the technique adopted in~\citep{lian2015asynchronous} for SGD. Furthermore, as an asynchronous algorithm, AsyVRSC-Shared incurs delays inevitably. When one worker has read the parameter and is computing the gradient, other workers may have finished their computation and updated the parameter in the shared memory. Therefore, the parameter one worker reads from the shared memory is delayed. We use $x_{t-\tau_t^s}^s$ to denote the delayed parameter used for computing the gradient of the \emph{t-th} inner update in the \emph{s-th} outer loop, where $\tau_t^s$ indicates the time delay.

As described in Algorithm~\ref{alg:shared}, in the \emph{t-th} inner loop, we uniformly sample a subscription $k_{t}$ from $[d_{1}]$ and update $x_{t}^{s}$ by:
{\small
\begin{equation}
(x_{t+1}^{s})_{k_{t}} = (x_{t}^{s})_{k_{t}} - \eta (\nabla \widehat{f}(x_{t-\tau_t^s}^{s}))_{k_{t}},
\end{equation}
}where $\eta$ denotes the learning rate, $\nabla \widehat{f}(x_{t-\tau_t^s}^{s})$ denotes the delayed variance reduced stochastic gradient used for the \emph{t-th} inner update in the \emph{s-th} outer loop, which is computed with (\ref{in_G}), (\ref{in_delta_G}) and (\ref{in_delta_f}). In addition, if the time delays have an upper bound $T$, we can represent $x_{t-\tau_t^s}^{s}$ as:
{\small
\begin{equation}
x_{t-\tau_t^s}^{s} = x_{t}^{s} - \sum_{j \in J(t)} (x_{j+1}^{s} - x_{j}^{s}),
\end{equation}
}where $J(t) \subseteq \{t,t-1,...,t-T + 1\}$ is a subset of previous iterations.

As $x_{r}^{s}$ in line 13 of Algorithm~\ref{alg:shared} is mainly used for theoretical analysis, we can replace it with $x_{K}^{s}$ in practice. This simplifies the computation and we have not observed much difference in the convergence speed empirically.
For efficient implementation, we can sample $k_t$ before line 9, then we only need to compute the corresponding part of 
$\nabla \widehat{f}(x_{t-\tau_t^s}^{s})$.

\subsection{AsyVRSC-Distributed}

\begin{algorithm}[!htb]
\caption {AsyVRSC-Distributed in the Master Node}
\label{alg:distributed_master}
\hspace*{0.02in}{\bf Input:}
Inner iteration number $K$, outer iteration number $S$, worker number $W$, learning rate $\eta$, initial point $\widetilde{x}^{1} \in \mathbb{R}^{d_{1}}$.

\begin{algorithmic}[1]
\For{$s=1,2,...,S$}
\State \textbf{Phase 1:}
    \State Broadcast $\widetilde{x}^s$ to all workers
    \State Receive and aggregate $G^{k}(\widetilde{x}^s)$ and $\nabla G^{k}(\widetilde{x}^s)$: 
    $$G(\widetilde{x}^s)= \frac{1}{n_2} \sum_{k=1}^{W} G^{k}(\widetilde{x}^s), \nabla G(\widetilde{x}^s) = \frac{1}{n_2} \sum_{k=1}^{W} \nabla G^{k}(\widetilde{x}^s)$$
    \State Broadcast $G(\widetilde{x}^s)$ to all workers
    \State Receive and aggregate $\nabla F^{k}(G(\widetilde{x}^s))$:
    $$\nabla F(G(\widetilde{x}^s)) = \frac{1}{n_1} \sum_{k=1}^{W}\nabla F^{k}(G(\widetilde{x}^s))$$
    \State Compute the full gradient $\nabla f(\widetilde{x}^s)$ using (\ref{full_grad})
    \State {\small  $x_{0}^{s} = \widetilde{x}^s$ } 
    \State Broadcast $\nabla G(\widetilde{x}^s), \nabla f(\widetilde{x}^s)$ and $x_{0}^{s}$ to all workers
\State \textbf{Phase 2:}
    \For{$t=0,1,2,...,K-1$}
        \State Receive gradient $\nabla \widehat{f}(x_{t-\tau_t^s}^{s})$ from one worker $w_t$
        \State Update $x_{t+1}^{s} = x_{t}^{s} - \eta \nabla \widehat{f}(x_{t-\tau_t^s}^{s})$
        \State Send $x_{t+1}^{s}$ to the worker $w_t$
    \EndFor
\State $\widetilde{x}^{s+1} = x_{r}^{s}$ for randomly chosen $r \in \{0,\cdots,K-1\}$
\EndFor
\State \Return $\widetilde{x}^{S+1}$
\end{algorithmic}
\end{algorithm}

\begin{algorithm}[!htb]
\caption {AsyVRSC-Distributed in the \emph{k-th} Worker Node}
\label{alg:distributed_worker}
\hspace*{0.02in}{\bf Input:}
Mini-batch size $a, b$.

\begin{algorithmic}[1]
\For{$s=1,2,...,S$}
\State \textbf{Phase 1:}
    \State Receive $\widetilde{x}^s$ from the master
    \State Compute $G^{k}(\widetilde{x}^{s})$ and $\nabla G^{k}(\widetilde{x}^{s})$ using (\ref{k-th_G}) and (\ref{k-th_delta_G}) and send them to the master.
    \State Receive $G(\widetilde{x}^s)$ from the master
    \State Compute $\nabla F^{k}(G(\widetilde{x}^{s}))$ using (\ref{k-th_delta_F}) and send it to the master
    \State Receive $\nabla G(\widetilde{x}^s),\nabla f(\widetilde{x}^s)$ from the master
\State \textbf{Phase 2:}
    \State Receive $x_{t- \tau_{t}^{s}}^{s}$ from the master
    \State Uniformly sample $A_t$ and $B_t$ from $[n_2]$ with replacement, where $|A_t|=a, |B_t|=b$
    \State Compute $\nabla \widehat{f}(x_{t-\tau_t^s}^{s})$ using (\ref{in_G}), (\ref{in_delta_G}) and (\ref{in_delta_f})
    \State Send $\nabla \widehat{f}(x_{t-\tau_t^s}^{s})$ to the master
\EndFor
\end{algorithmic}
\end{algorithm}

In a master-worker architecture, suppose that there are a master node and $W$ worker nodes. The master maintains the parameter and updates it when receiving a gradient from any worker. Each local worker pulls the current parameter from the master, calculates the gradient locally and sends it to the master independently~\citep{agarwal2011distributed}. 

In AsyVRSC-Distributed, similar to AsyVRSC-Shared, there are two phases in each outer loop. The description of AsyVRSC-Distributed is presented in Algorithm~\ref{alg:distributed_master} and Algorithm~\ref{alg:distributed_worker}, which show the operations of the master node and the \emph{k-th} worker node, respectively.

\textbf{Phase 1:} As shown in line $3$-$9$ in Algorithm~\ref{alg:distributed_master} and line $3$-$7$ in Algorithm~\ref{alg:distributed_worker}, at the beginning of the \emph{s-th} outer loop, the master broadcasts $\widetilde{x}^s$ to all workers, then the workers calculate the full gradient collectively. Specifically, we equally divide $[n_1]$ and $[n_2]$ into $W$ blocks and use $N_k$ and $M_k$ to denote the \emph{k-th} blocks, respectively.
The \emph{k-th} worker calculates the corresponding parts belonging to $N_k$ and $M_k$:
{\small \begin{equation}\label{k-th_G}
G^{k}(\widetilde{x}^s) = \sum_{j \in N_k} G_{j} (\widetilde{x}^s),
\end{equation}
\begin{equation}\label{k-th_delta_G}
\nabla G^{k}(\widetilde{x}^s) = \sum_{j \in N_k} \nabla G_{j} (\widetilde{x}^s),
\end{equation}
\begin{equation}\label{k-th_delta_F}
\nabla F^{k}(G(\widetilde{x}^s)) = \sum_{i \in M_k} \nabla F_{i} (G(\widetilde{x}^s)),
\end{equation} 
}After that, the master aggregates the gradients from all workers to get the full gradient and broadcasts $G(\widetilde{x}^s),\nabla G(\widetilde{x}^s)$ and $\nabla f(\widetilde{x}^s)$ to all workers. 

\textbf{Phase 2:} As shown in line $11$-$15$ in Algorithm~\ref{alg:distributed_master} and line $9$-$12$ in Algorithm~\ref{alg:distributed_worker}, in the inner loops, all workers compute the variance reduced stochastic compositional gradient in an asynchronous way and the master conducts the updates.

In a master-worker architecture, the updates are atomic 
if we let the master only respond to a single worker in each iteration. 
The atomic operation ensures that the parameter one worker gets from the master is a real state of $x$.
Similar to AsyVRSC-Shared, the parameter used to compute the gradient in AsyVRSC-Distributed may be delayed. We also use $\tau_{t}^{s}$ to denote the time delay. As described in Algorithm~\ref{alg:distributed_master}, when receiving a gradient from a worker, the master updates $x_{t}^{s}$ by:
\begin{equation}
 x_{t+1}^{s} = x_{t}^{s} - \eta \nabla \widehat{f}(x_{t-\tau_t^s}^{s}),
\end{equation} 
where $\nabla \widehat{f}(x_{t-\tau_t^s}^{s})$ is computed with (\ref{in_G}), (\ref{in_delta_G}) and (\ref{in_delta_f}).

Similar to AsyVRSC-Shared, $x_{r}^{s}$ in line $16$ of Algorithm~\ref{alg:distributed_master} can be replaced with $x_{K}^{s}$ in practice.

\section{Theoretical Analysis}
\label{theory}
In this section, we give the theoretical analysis for the two algorithms proposed in the previous section. The main difficulties in the theoretical analysis of asynchronous algorithms are caused by the time delays. We prove that AsyVRSC-Shared and AsyVRSC-Distributed can achieve linear speedup when the time delays can be bounded by the data dimensionality and the sparsity ratio of the partial gradients, respectively.
\subsection{AsyVRSC-Shared}
At first, we introduce some basic assumptions, which are commonly used in theoretical analysis for composition optimization~\citep{wang2017stochastic,lian2017finite,huo2017accelerated,yu2017fast}:
\newtheorem{assu}{Assumption}
\begin{assu}[Lipschitz Gradient] \label{assumption:smoothness} 
There exist Lipschitz constants $L_F$, $L_G$ and $L_f$ for $ \nabla F(x),\nabla G(x)$ and $\nabla f(x)$, respectively, such that for $ \forall i \in [n_1], \forall j \in [n_2] $, $\forall x_1, x_2 \in \mathbb{R}^{d_1}$ and $\forall y_1, y_2 \in \mathbb{R}^{d_2}$:
{\small
\begin{equation} \label{assum1:1}
\|\nabla F_{i}(y_1) - \nabla F_{i}(y_2)\| \leq L_{F}\|y_1-y_2\|,
\end{equation}
\begin{equation} \label{assum1:2}
\|\nabla G_{j}(x_1)-\nabla G_{j}(x_2)\| \leq L_{G} \|x_1-x_2\|,
\end{equation}
\begin{equation} \label{assum1:3}
\|\nabla f_{i,j}(x_1) - \nabla f_{i,j}(x_2)\| \leq L_{f} \|x_1-x_2\|,
\end{equation} 
}where $\nabla f_{i,j}(x) = (\nabla G_{j}(x))^{T} \nabla F_{i}(G(x)) $.
From (\ref{assum1:3}) we immediately have:
{\small
\begin{equation} \label{assum1:4}
\|\nabla f(x_1) - \ \nabla f(x_2)\| \leq L_{f} \|x_1-x_2\|.
\end{equation}
}
\end{assu}

\newtheorem{strongly}[assu]{Assumption}
\begin{strongly}[Strong Convexity] \label{assumption:strongly convex}
$f(x)$ is a strongly convex function with parameter $\mu_{f}$. For $\forall x_1,x_2 \in \mathbb{R}^{d_1}$:
{\small
\begin{equation}
f(x_1) - f(x_2) \geq \langle \nabla f(x_2), x_1-x_2\rangle + \frac{\mu_{f}}{2} \|x_1 - x_2\|^{2}.
\end{equation}
}Furthermore, if $f(x)$ is strongly convex, there exists an unique optimal solution $x^{*}$ to problem (\ref{basic_formula}).
\end{strongly}

\newtheorem{bounded}[assu]{Assumption}

\begin{bounded}[Bounded Gradient] \label{assumption:bounded gradient}
The gradients $\nabla F_{i}(x)$ and $\nabla G_{j}(x)$ are bounded by constants $B_{F}$ and $B_{G}$ respectively. For $\forall i \in [n_1]$, $\forall j \in [n_2]$, $\forall x_1, x_2 \in \mathbb{R}^{d_1}$ and $\forall y_1, y_2 \in \mathbb{R}^{d_2}$:
{\small
\begin{eqnarray}
\|\nabla F_{i}(y_1)\| &\leq& B_{F},
\\
\|\nabla G_{j}(x_1)\| &\leq& B_{G},
\end{eqnarray}
}then $F_{i}(x)$ and $G_{j}(x)$ are Lipschitz functions that satisfy:
{\small
\begin{eqnarray}
\| F_{i}(y_1) - F_{i}(y_2)\| &\leq& B_{F} \|y_1-y_2\|,
\\
\|G_{j}(x_1) - G_{j}(x_2)\| &\leq& B_{G} \|x_1-x_2\|,
\end{eqnarray}
}
\end{bounded}

In asynchronous parallel algorithms, the gradients used for updating may be delayed. It is natural to assume an upper bound for the time delays:
\newtheorem{bounded_delay}[assu]{Assumption}
\begin{bounded_delay}[Bounded Delay] \label{assumption:bounded delay} 
Assume that there exits a constant $T$ such that $\tau_{t}^{s} \leq T$, for any outer loop $s$ and inner loop $t$. In practice, $T$ is roughly proportional to the number of workers.
\end{bounded_delay}

We first propose two lemmas to bound the variance of the gradients and the variance of the estimated inner function values. They are the cornerstones of all subsequent analysis. 

\begin{lemma} \label{lemma:1}
Let $x^{*}$ be the optimum to problem (\ref{basic_formula}) such that $x^{*} = arg min_{x \in \mathbb{R}^{d_{1}}} f(x)$. Under Assumptions 1-3, the following inequality holds:
{\small 
\begin{equation} \label{lemma:bound_variance} 
\mathbb{E}\|\nabla \widehat{f}(x_{t}^{s}) \|^2 \leq R \mathbb{E}(f(x_{t}^{s}) - f(x^{*})+ f(\widetilde{x}^{s}) - f(x^{*})),
\end{equation} 
}where $R=\frac{64}{\mu_{f}}(\frac{B_{G}^4 L_{F}^2}{a} + \frac{B_{F}^{2} L_{G}^{2}}{b} ) + 8 L_{f}$.
\end{lemma}

\begin{lemma} \label{lemma:2}
Under Assumptions 1-3, the variance of the estimated inner function values can be bounded as following:
{\small
\begin{equation} \label{lemma:bound_G}
\mathbb{E} \|\widehat{G}_{t - \tau_{t}^{s}}^{s} - G(x_{t - \tau_{t}^{s}}^{s})\|^2 \leq \frac{3B_{G}^2}{a} \mathbb{E} (\|\widetilde{x}^{s} - x^{*}\|^2 + \|x_{t}^{s} - x^{*}\|^2 + \|x_{t}^{s}- x_{t - \tau_{t}^{s}}^{s}\|^2)
\end{equation}
}
\end{lemma}

As shown in Lemma~\ref{lemma:1} and Lemma~\ref{lemma:2}, $\mathbb{E}\|\nabla \widehat{f}(x_{t}^{s})\|^2$ can be bounded by the optimality gap of $f(x)$, defined by $f(x)-f(x^*)$, and the variance of the estimated inner function values can be bounded by the Euclidean distance from the current parameter to the optimum. This means the variance of AsyVRSC asymptotically goes to zero as $x_t^s$ and $\widetilde{x}^s$ converge to $x^*$, and it is the main reason why AsyVRSC can converge with a constant learning rate.

One of the main difficulties in the analysis of asynchronous algorithms is to bound the delayed gradients. We derive the upper bound of the delayed gradients in the following lemmas.

\begin{lemma} \label{lemma:3}  
Assume Assumptions $1$-$3$ hold. The delayed estimated gradient of the inner function can be bounded as following:
{\small
\begin{equation}
\mathbb{E}\|\nabla \widehat{G}_{t - \tau_{t}^{s}}^{s}\|^2 \leq 10 B_{G}^2.
\end{equation}
}
\end{lemma}

\begin{lemma} \label{lemma:4}
Assume Assumptions $1$-$4$ hold. In each epoch of AsyVRSC-Shared, the sum of all delayed gradients can be bounded by the sum of undelayed gradients as following:
{\small
\begin{equation} \label{lemma:bound_delayed_shared}
\sum_{t=0}^{K-1} \mathbb{E} \|\nabla \widehat{f}(x_{t - \tau_{t}^{s}}^{s})\|^2 \leq \frac{2}{1- [40B_{G}^4 L_{F}^2 + 4 B_{F}^2 L_{G}^2] \frac{T^2 \eta^2}{d_1} } \sum_{t=0}^{K-1} \mathbb{E} \|\nabla \widehat{f}(x_{t}^{s})\|^2.
\end{equation}
}
\end{lemma}

Combining these lemmas, we obtain the following theorem.

\newtheorem{theorem_shared_1}{Theorem}
\begin{theorem_shared_1}  \label{theorem_shared_1:1}
Under Assumptions $1$-$4$, AsyVRSC-Shared has geometric convergence in expectation:
{\small
\begin{equation}\label{theorem_shared_1:2}
\mathbb{E}[f(\widetilde{x}^{s+1}) - f(x^{*})] \leq \frac{ \frac{2}{\mu_f} + PQRK + U}{ \frac{7\eta K}{4d_1} - PQRK - U} \mathbb{E}[f(\widetilde{x}^{s}) - f(x^{*})],
\end{equation} 
}where
{\small
\begin{eqnarray} \label{theorem_shared_1:parameter}
U &=& \frac{48 \eta B_{G}^{4} L_{F}^{2} K}{d_1 a \mu_{f}^2}, 
\nonumber\\
P &=&\frac{\eta^2}{d_1} + \frac{L_{f} T^2 \eta^3}{d_{1}^2} + \frac{24 \eta^3 B_{G}^4 L_{F}^2 T^2}{d_{1}^2 a \mu_{f}}, 
\nonumber\\
Q &=&  \frac{2}{1-[40B_{G}^4L_{F}^2 + 4B_{F}^2L_{G}^2]\frac{\eta^2 T^2}{d_1}},
\nonumber\\
R &=& \frac{64}{\mu_{f}}(\frac{B_{G}^4 L_{F}^2}{a} + \frac{B_{F}^2 L_{G}^{2}}{b})+ 8 L_{f}.
\end{eqnarray}
}
\end{theorem_shared_1}

\newtheorem{corollary_shared_1}{Corollary}
\begin{corollary_shared_1} \label{corollary_shared_1:1}
Suppose that the conditions in Theorem 1 hold and we set the parameters as $a= \max\{ \frac{1024 B_{G}^4 L_{F}^2}{\mu_{f}^2}, \frac{32B_{G}^4 L_F^2}{5\mu_f L_f}\}$, $b = \frac{32 B_{F}^2 L_{G}^2}{\mu_{f} L_{f}}$, $K = \frac{1024 L_{f} d_1}{\mu_{f}}$, and $\eta = \min \{ \frac{1}{9 B_{G}^2 L_{F}}, \frac{1}{9B_{F} L_{G}},\frac{1}{320 L_{f}}\}$. Then if $T$ can be bounded:
{\small
$$T \leq \sqrt{d_1},$$
}AsyVRSC-Shared has the following linear convergence rate:
{\small
\begin{equation}
\mathbb{E}[f(\widetilde{x}^{s+1}) - f(x^{*})] \leq \frac{2}{3} \mathbb{E}[f(\widetilde{x}^{s}) - f(x^{*})].
\end{equation} 
}
\end{corollary_shared_1}

From Corollary 1, if we choose $a$, $b$, $K$ and $\eta$ properly and $T$ can be bounded by $\sqrt{d_1}$, AsyVRSC-Shared has a linear convergence rate. Therefore, if we want to achieve $\mathbb{E} f(\widetilde{x}^s) - f(x^*) \leq \epsilon$, the number of updates we need to take is $O((n_1+n_2+K(a+b))\log{\frac{1}{\epsilon}}) = O((n_1+n_2+\kappa^3d_1)\log{\frac{1}{\epsilon}})$, where $\kappa = \max \{ \frac{L_F}{\mu_f}, \frac{L_G}{\mu_f}, \frac{L_f}{\mu_f} \}$ denotes the condition number of the loss function. Since we only calculate one single component of the gradient at each inner iteration, the overall query complexity is $O((n_1+n_2+\kappa^3) \log{\frac{1}{\epsilon}})$, which is independent with the number of workers. The query complexity is consistent with the theoretical result of VRSC~\citep{lian2017finite}. Since the number of updates we need to take in the parallel AsyVRSC-Shared is the same as that in the single-machine VRSC in order and the constant is irrelevant to the number of workers, it follows that AsyVRSC enjoys linear speedup of parallel computation.

\subsection{AsyVRSC-Distributed}
In this subsection, we give a theoretical analysis of AsyVRSC-Distributed. To ensure the linear convergence rate of AsyVRSC-Distributed, we need to make a further sparsity assumption:

\newtheorem{sparsity}[assu]{Assumption}
\begin{sparsity}[Sparsity] \label{assumption:sparsity}

$\nabla F_{i}(x)$, $\nabla G_{j}(x)$ and $\nabla f_{ij}(x)$  are all sparse, where $\nabla f_{ij}(x) = (\nabla G_{j}(x))^{T} \nabla F_{i}(x)$. We introduce  $\Delta_{F} = \max_{i \in [n_1]} | \{k \mid \nabla_{k} F_{i}(x) \neq 0, k \in [d_2]\}$, $\Delta_{G} = \max_{j \in [n_2]} | \{ (k,l) \mid \nabla_{k,l} G_{j}(x) \neq 0, k \in [d_2], l \in [d_1]\} |$ and $\Delta_{f} = \max_{i \in [n_1], j \in [n_2]} | \{k \mid \nabla_k f_{ij}(x) \neq 0, k \in [d_1]\} |$, where $1 \leq \Delta_F \leq d_2$, $1 \leq \Delta_G \leq d_1 d_2$ and $1 \leq \Delta_f \leq d_1$. Then we define $\Delta = \max\{\frac{\Delta_F }{d_2}, \frac{\Delta_G}{d_1 d_2}, \frac{\Delta_f}{d_1}\}$. We always have $ \frac{1}{d_1 d_2} \leq \Delta \leq 1$.
\end{sparsity}

The sparsity assumption is common in analyzing distributed variance reduced asynchronous algorithms to solve strongly convex problems; see previous works \citep{reddi2015variance,meng2017asynchronous,leblond2017asaga,pedregosa2017breaking}. AsyVRSC-Distributed also features in variance reduction, and hence inherits this assumption.
Roughly speaking, in the analysis of variance reduced asynchronous algorithms, the time delay appears as a dominating factor in the rate of convergence. To obtain a favorable rate, one has to handle the time delay with assumptions such as bounded delay and sparsity.

Indeed, the sparsity assumption is satisfied in many applications. For example, both portfolio management problem~\citep{lian2017finite} and on-policy learning problem~\citep{wang2016accelerating} satisfy this assumption when the dataset is sparse.

Similar to the theoretical analysis in Section $4.1$, we derive a lemma to bound the delayed gradients.

\begin{lemma}\label{lemma:5}
Assume Assumptions $1$-$5$ hold. In each epoch of AsyVRSC-Distributed, the sum of all delayed gradients can be bounded by the sum of undelayed gradients as following:
{\small
\begin{equation} \label{lemma:bound_delayed_distributed}
\sum_{t=0}^{K-1} \mathbb{E} \|\nabla \widehat{f}(x_{t - \tau_{t}^{s}}^{s})\|^2 \leq \frac{2}{1- [40B_{G}^4 L_{F}^2 + 4 B_{F}^2 L_{G}^2] \Delta \eta^2 T^2} \sum_{t=0}^{K-1} \mathbb{E} \|\nabla \widehat{f}(x_{t}^{s})\|^2.
\end{equation}
}
\end{lemma}

\newtheorem{theorem_distribute_1}[theorem_shared_1]{Theorem}
\begin{theorem_distribute_1}\label{theorem_ditribute_1:1}
Under Assumptions 1-5, AsyVRSC-Distributed has geometric convergence in expectation:
{\small
\begin{equation} \label{theorem_distribution_1:2}
\mathbb{E}[f(\widetilde{x}^{s+1}) - f(x^{*})] \leq \frac{\frac{\mu_f}{2} + PQRK +U}{\frac{7}{4}\eta K - PQRK - U} \mathbb{E}[f(\widetilde{x}^{s}) - f(x^{*})],
\end{equation} 
}where
{\small
\begin{eqnarray} \label{theorem_distribution_1:parameter}
U &=& \frac{48 \eta B_{G}^{4} L_{F}^{2} K}{a \mu_{f}^2},
\nonumber\\
P &=&\eta^2 + \Delta L_{f} T^2 \eta^3 + \frac{24 \eta^3 B_{G}^4 L_{F}^2 \Delta T^2}{a \mu_{f}},
\nonumber\\
Q &=&  \frac{2}{1-[40B_{G}^4L_{F}^2 + 4B_{F}^2L_{G}^2]\Delta \eta^2 T^2},
\nonumber\\
R &=& \frac{64}{\mu_{f}}(\frac{B_{G}^4 L_{F}^2}{a} + \frac{B_{F}^2 L_{G}^{2}}{b})+ 8 L_{f}.
\end{eqnarray} 
}
\end{theorem_distribute_1}

\newtheorem{corollary_distribute_1}[corollary_shared_1]{Corollary}

\begin{corollary_distribute_1} \label{corollary_ditribute_1:1}
Suppose the conditions in Theorem 2 hold and we set the parameters as $a = \max\{\frac{1024 B_{G}^4 L_{F}^2}{\mu_{f}^2}, \frac{32B_G^4L_F^2}{5\mu_f L_f}\}$, $b = \frac{32 B_{F}^2 L_{G}^2}{\mu_{f} L_{f}}$, $K = \frac{1024 L_{f}}{\mu_{f}}$ and $\eta = \min \{ \frac{1}{9 B_{G}^2 L_{F}}, \frac{1}{9B_{F} L_{G}},\frac{1}{320 L_{f}}\}$.
Then if $T$ can be bounded:
{\small
$$T \leq \sqrt{\Delta^{-1}},$$ 
}AsyVRSC-Distributed has the following linear convergence rate:
{\small
\begin{equation}
\mathbb{E}[f(\widetilde{x}^{s+1}) - f(x^{*})] \leq \frac{2}{3} \mathbb{E}[f(\widetilde{x}^{s}) - f(x^{*})].
\end{equation} 
}
\end{corollary_distribute_1}

Theorem 2 and Corollary 2 give a linear convergence analysis of AsyVRSC-Distributed, which depends on the sparsity ratio of $\nabla F_i (x)$, $\nabla G_j (x)$ and $\nabla f_{ij}(x)$. To achieve $f(x) - f(x^*) \leq \epsilon$, the number of updates we need to take is $O((n_1+n_2+K(a+b))\log{\frac{1}{\epsilon}})=O((n_1+n_2+\kappa^3) \log{\frac{1}{\epsilon}})$,  where $\kappa = \max \{ \frac{L_F}{\mu_f}, \frac{L_G}{\mu_f}, \frac{L_f}{\mu_f}\}$ denotes the condition number. Similar to AsyVRSC-Shared, AsyVRSC-Distributed can achieve linear speedup when $\Delta$ is small and $T \leq \sqrt{\Delta^{-1}}$.

\section{Experiments}
\label{exp}
In this section, we conduct experiments to verify the effectiveness of AsyVRSC-Shared and AsyVRSC-Distributed, including examples in reinforcement learning following~\citep{wang2016accelerating} and portfolio management following~\citep{lian2017finite}. For the shared-memory architecture, we use the OpenMP library\footnote{https://openmp.org/} to parallelize the multiple threads. As for the master-worker architecture, we use multi-process to simulate multi-machine operations. In real-world multi-machine operations, the speed and the speedup may be a little worse than that in our experiments due to the higher communication cost. The communications between the master and the workers are handled by the MPICH library\footnote{https://www.mpich.org/}. All the experiments are conducted on one single machine with 2 sockets, and each socket has 12 cores. Performance is evaluated by iteration speedup and running time speedup, which are defined as~\citep{lian2015asynchronous}:
{\small $$\textrm{Running~Time~Speedup} = \frac{\textrm{Running~time~of~using~one~worker}}{\textrm{Running~time~of~using~\emph{W}~workers}}.$$ }

{\small $$\textrm{Iteration~Speedup} = \frac{\textrm{Number~of~total~iterations~using~one~worker}}{\textrm{Number~of~total~iterations~using~\emph{W}~workers}} \times W.$$}

The thread number of AsyVRSC-Shared and the worker number of AsyVRSC-Distributed are both varying from 1 to 16. AsyVRSC-Shared with 1 thread and AsyVRSC-Distributed with 1 worker are approximately equivalent to VRSC~\citep{lian2017finite}. For AsyVRSC-Shared, updating only a single component of $x$ at each update is time-consuming in practice, therefore we uniformly select a subset of $x$ to conduct the update and we set the subset size to 10.
We tune the learning rate $\eta$ to get the best performance and empirically set the mini-batch sizes $a=b=5$.

\begin{figure}[!hbt]

\centering

\subfigure[Dataset-ref-1]{
\begin{minipage}[t]{0.3\linewidth}
\centering
\includegraphics[width=1.1\textwidth] {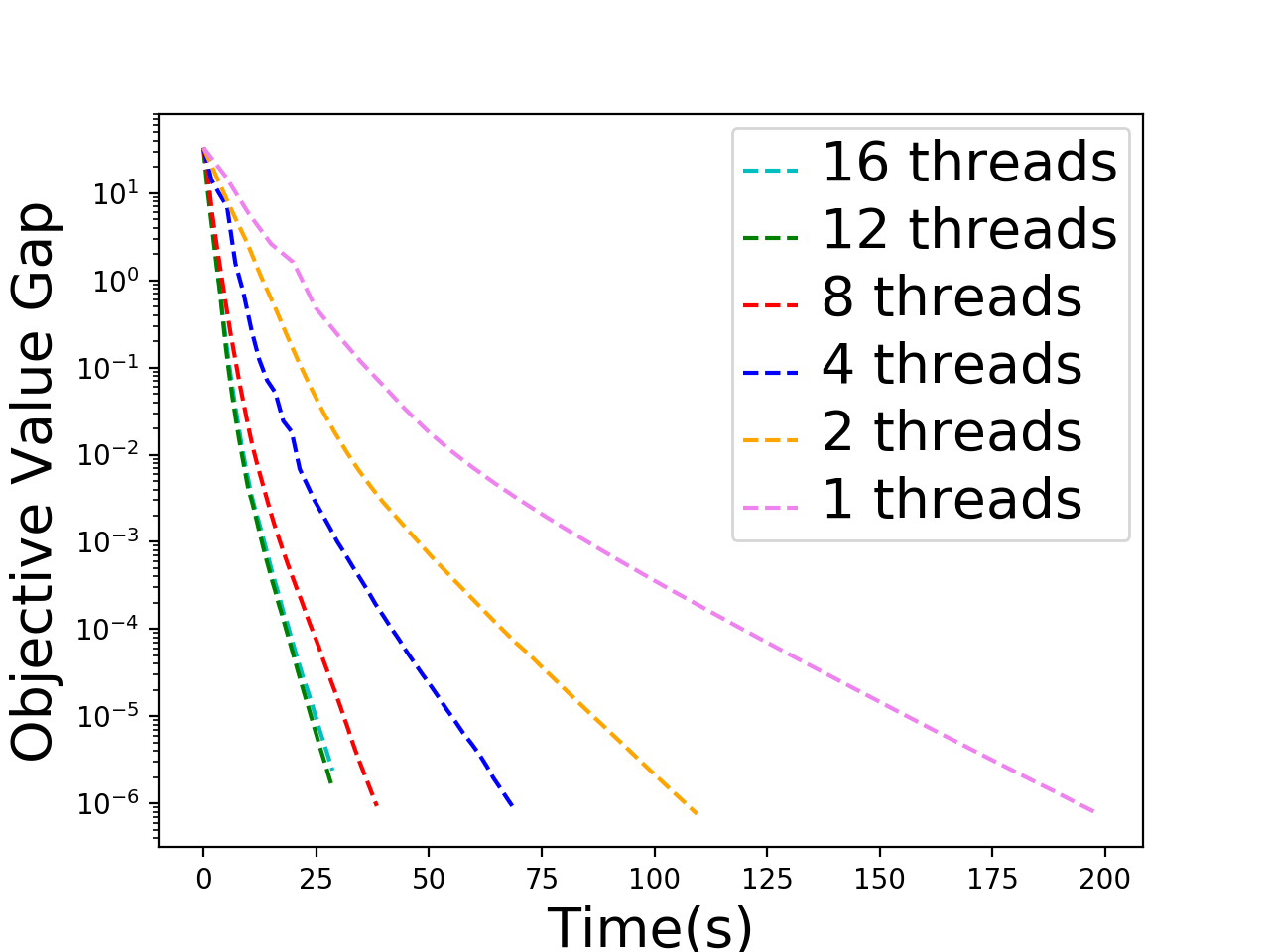}
\end{minipage}
}
\subfigure[Dataset-ref-1]{
\begin{minipage}[t]{0.3\linewidth}
\centering
\includegraphics[width=1.1\textwidth] {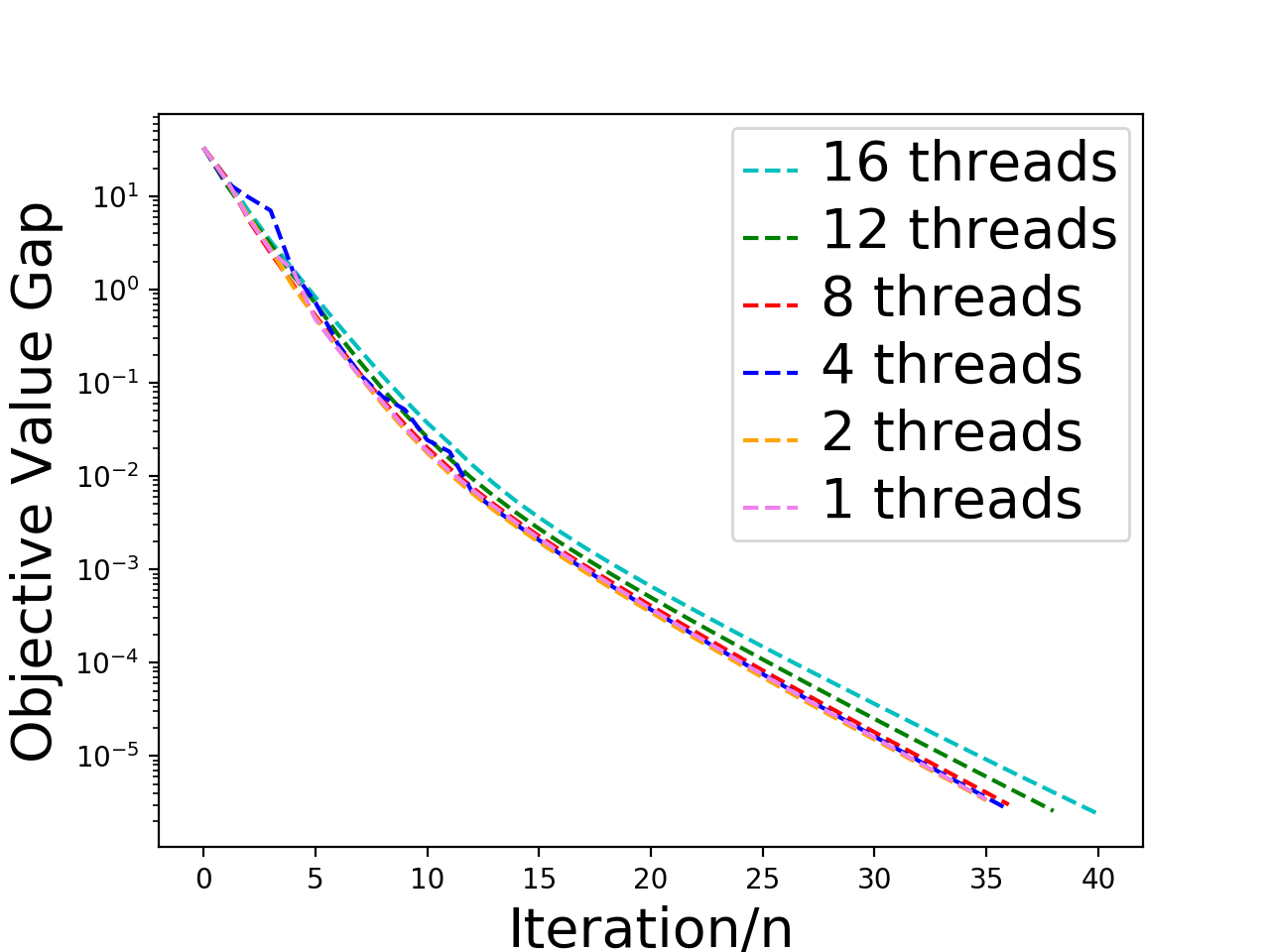}
\end{minipage}
}
\subfigure[Dataset-ref-1]{
\begin{minipage}[t]{0.3\linewidth}
\centering
\includegraphics[width=1.1\textwidth] {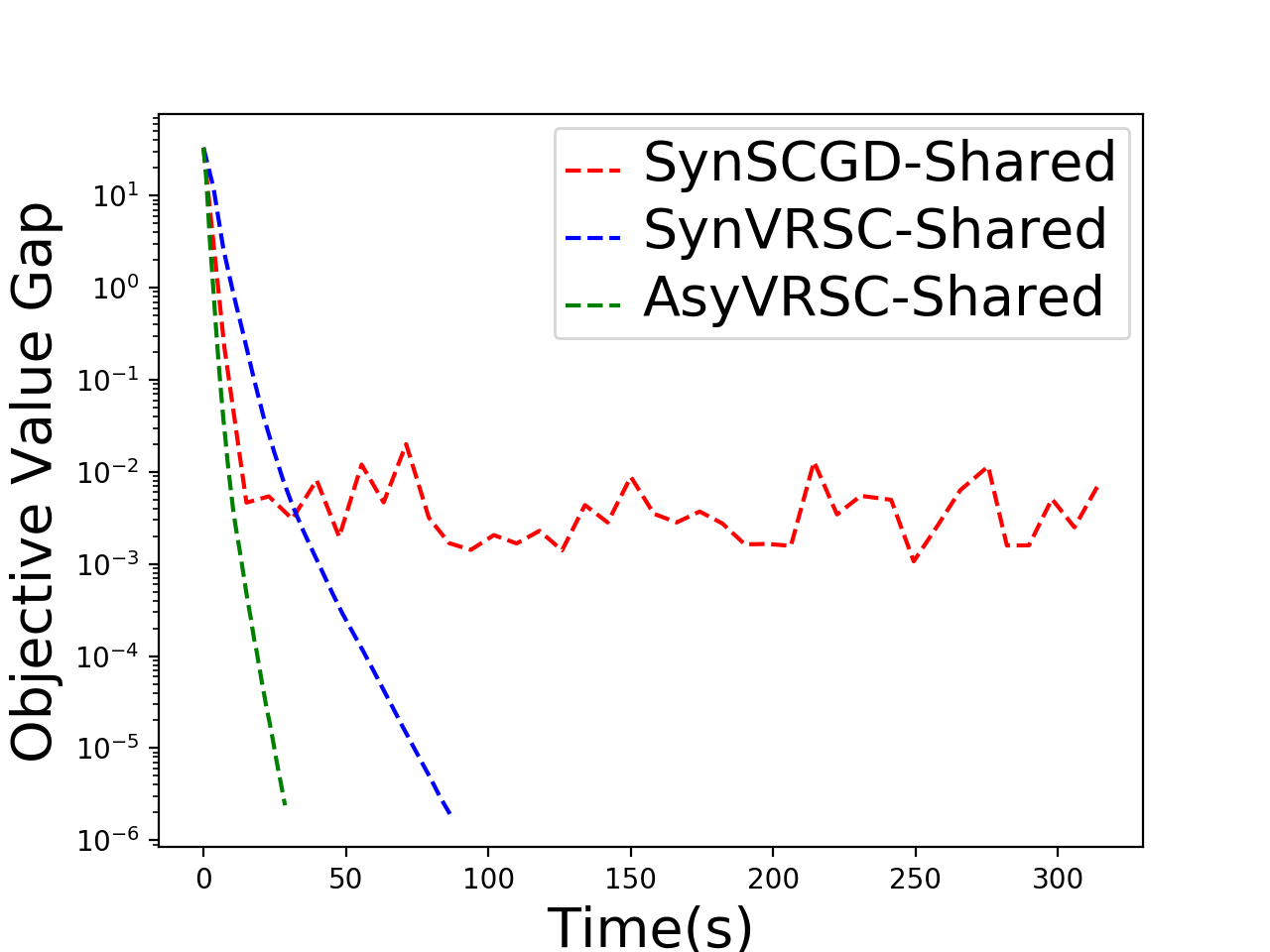}
\end{minipage}
}

\subfigure[Dataset-ref-2]{
\begin{minipage}[t]{0.3\linewidth}
\centering
\includegraphics[width=1.1\textwidth] {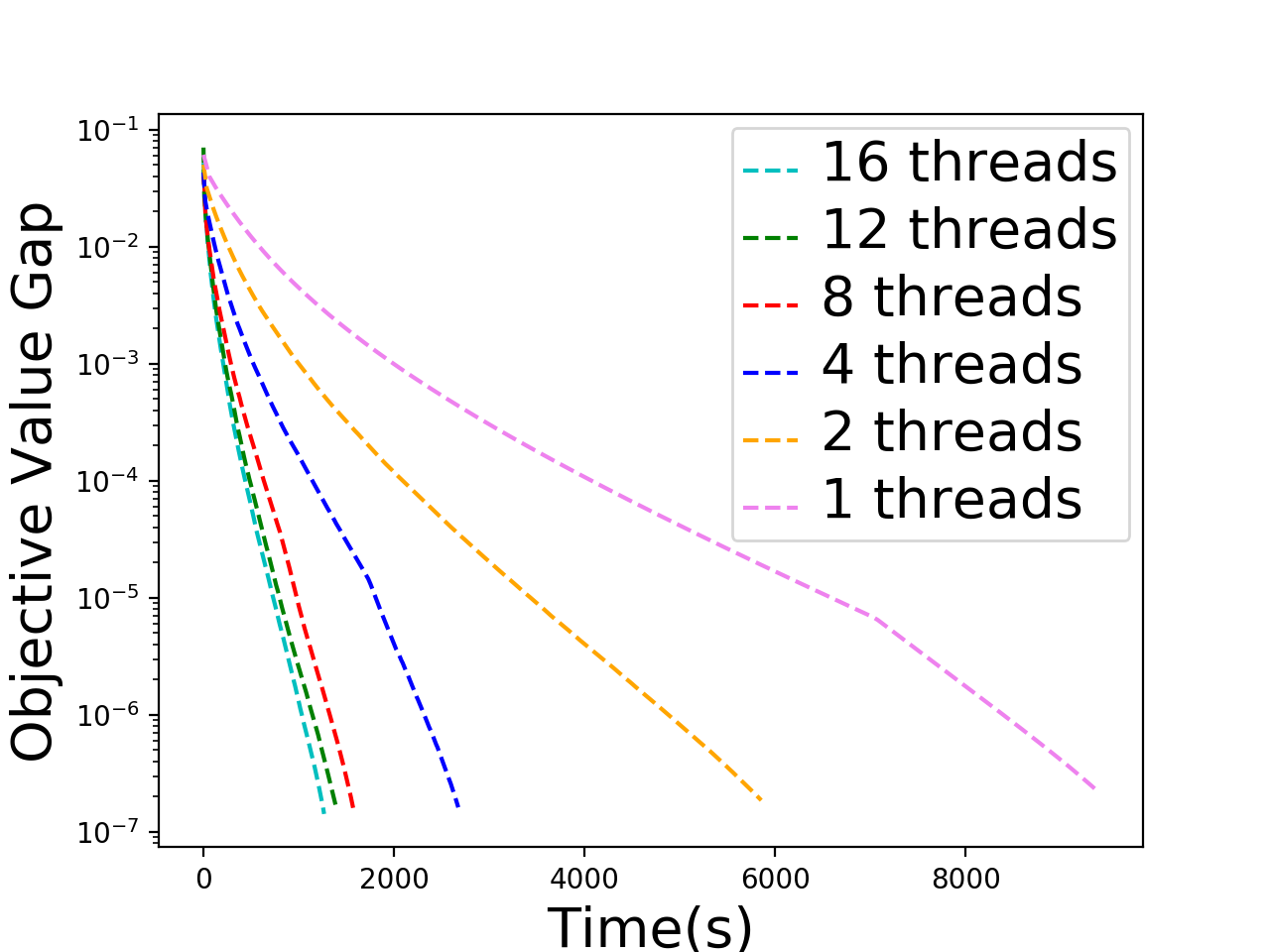}
\end{minipage}
}
\subfigure[Dataset-ref-2]{
\begin{minipage}[t]{0.3\linewidth}
\centering
\includegraphics[width=1.1\textwidth] {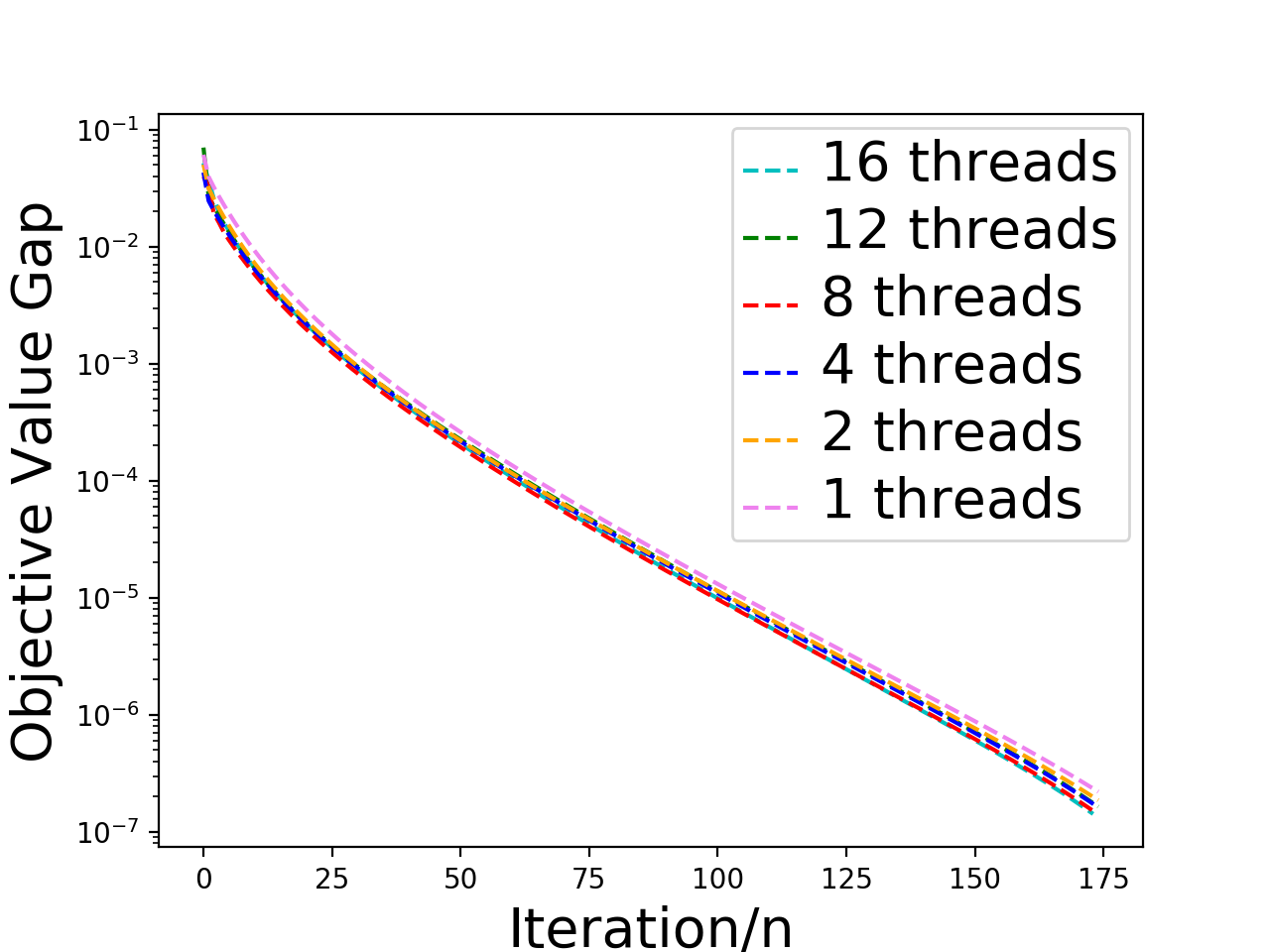}
\end{minipage}
}
\subfigure[Dataset-ref-2]{
\begin{minipage}[t]{0.3\linewidth}
\centering
\includegraphics[width=1.1\textwidth] {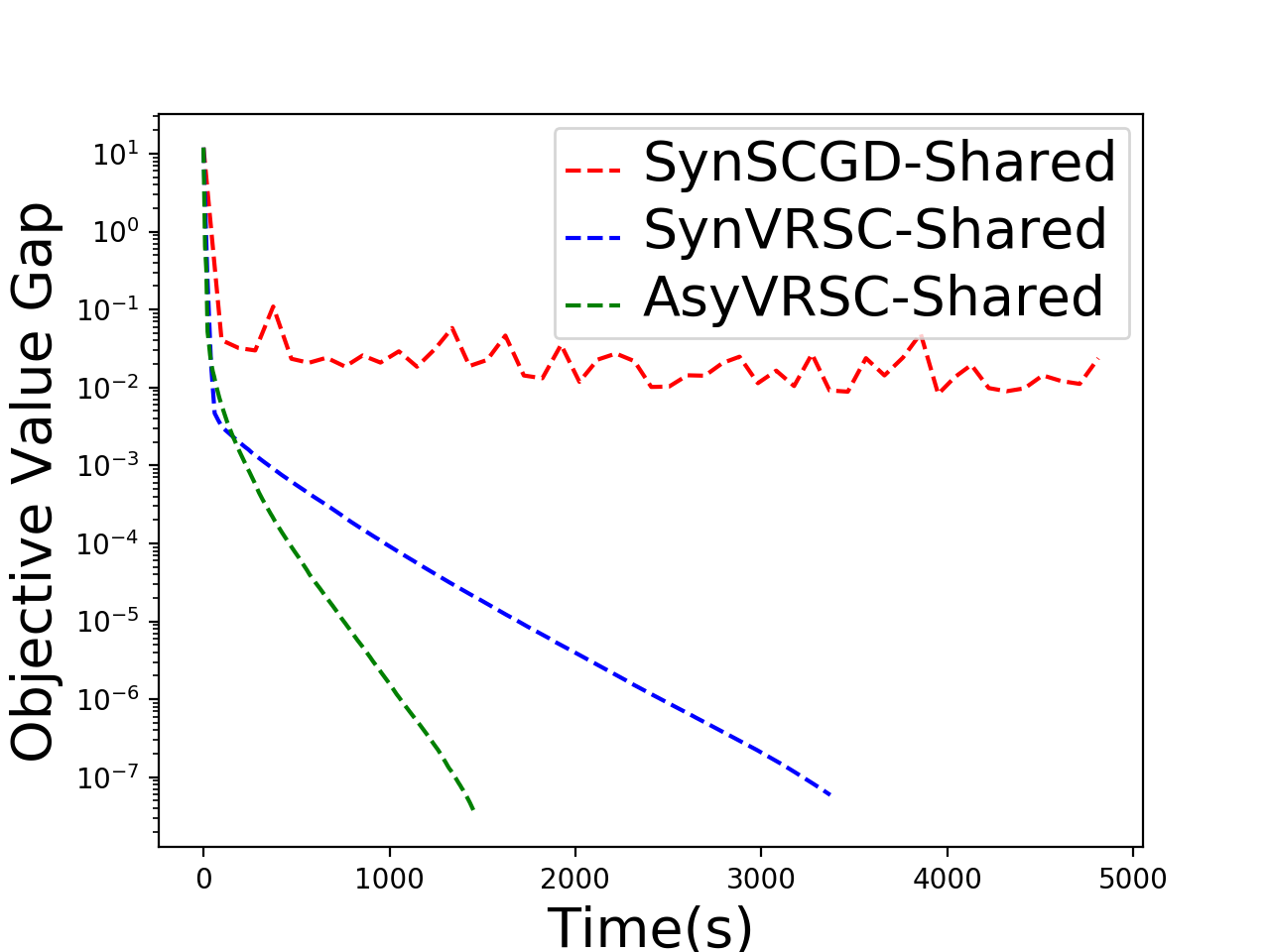}
\end{minipage}
}

\subfigure[Time Speedup]{
\begin{minipage}[t]{0.3\linewidth}
\centering
\includegraphics[width=1.1\textwidth] {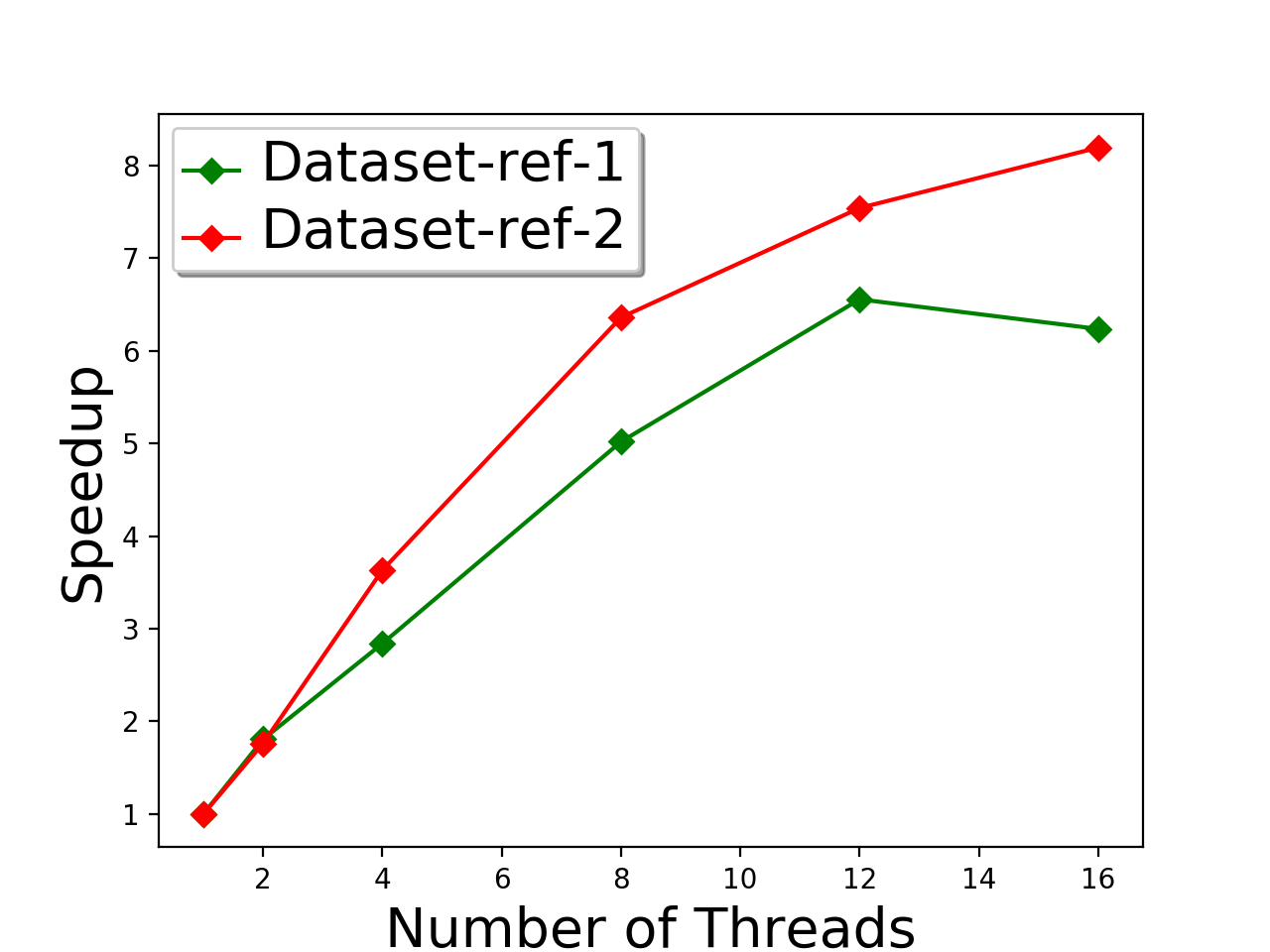}
\end{minipage}
}
\subfigure[Iteration Speedup]{
\begin{minipage}[t]{0.3\linewidth}
\centering
\includegraphics[width=1.1\textwidth] {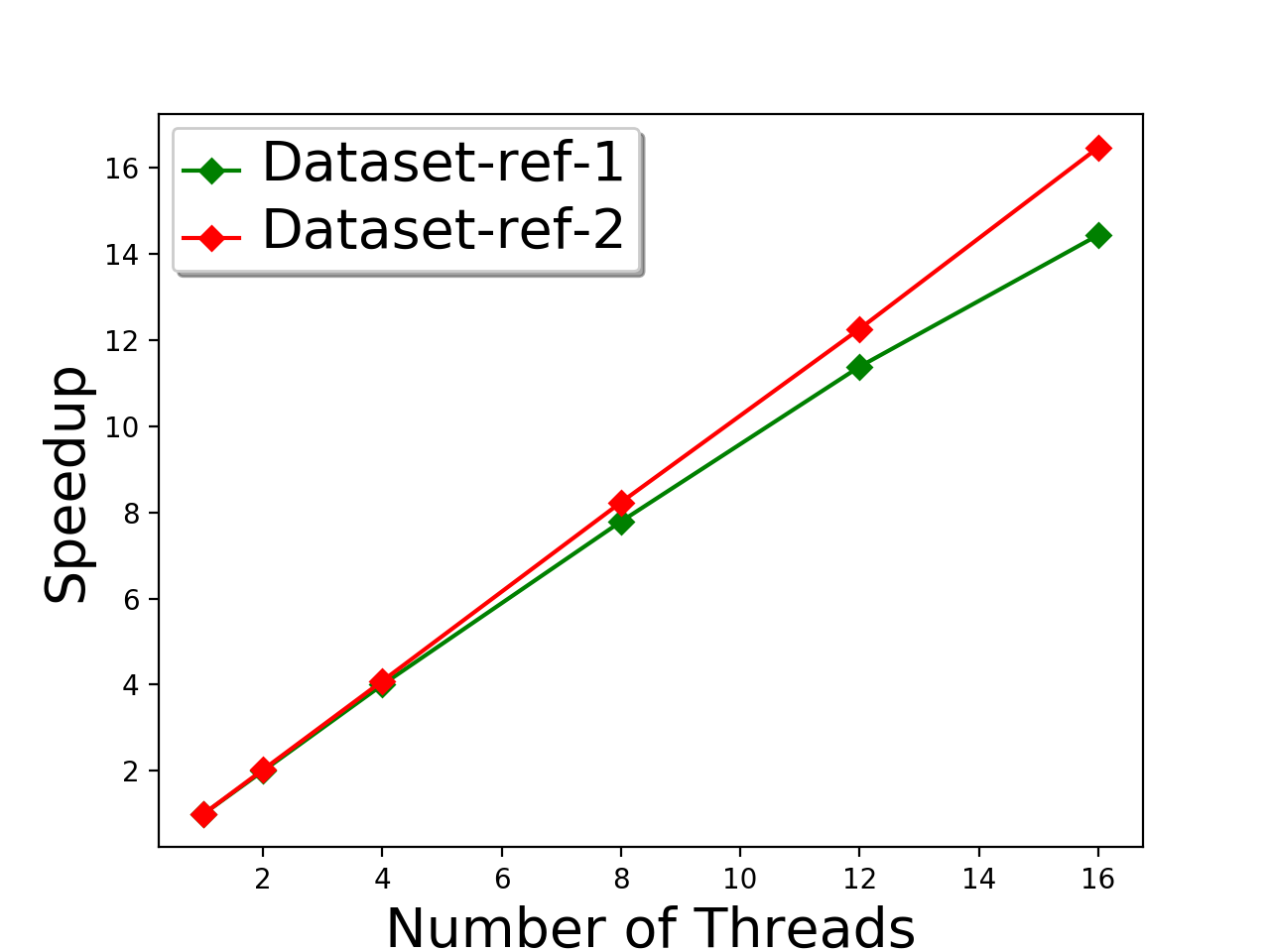}
\end{minipage}
}

\caption{Results for AsyVRSC-Shared on the reinforcement learning task. Figure (a)-(c) and (d)-(f) show the results for Dataset-ref-1 and Dataset-ref-2, respectively. Figure (c) and (f) show the comparision of AsyVRSC-Shared with SynSCGD-Shared and SynVRSC-Shared with 16 threads on the two datasets. Figure (g) and (h) show the Time Speedup and Iteration Speedup of AsyVRSC-Shared on the two datasets.}

\label{figure_ref}
\end{figure}

\subsection{Reinforcement Learning} \label{section_ref}
For on-policy learning, suppose that there are $S$ states and a fixed control policy $\pi$. The value function of each state can be approximated by an inner product of the state feature $\phi_s \in \mathbb{R}^d$ and target variable $x$, i.e., $V^{\pi}(s) = \phi_{s}^{T} x^*$. Then the on-policy problem can be formulated as:
{\small 
$$\min_{x \in \mathbb{R}^d} \frac{1}{S} \sum_{i=1}^{S}(\phi_{s}^{T} x - \sum_{s'}P_{s,s'}^{\pi}(r_{s,s'} + \gamma \cdot \phi_{s'}^{T} x))^2,$$
}where $\gamma \in (0,1)$ is a discount factor, $r_{s,s'}$ denotes the reward of transition from $s$ to $s'$ and $P_{s,s'}^{\pi}$ denotes the transition probability from state $s$ to state $s'$. We formulate this problem as a compositional problem in the form of (\ref{basic_formula})
by setting:
{\small
$$G_{j}(x) = (\phi_{1}^{T}x, S P_{1,j}^{\pi}(r_{1,j}+\gamma \phi_{j}^{T}x), \cdots, \phi_{S}^{T}x, S P_{S,j}^{\pi} (r_{S,j} + \gamma \phi_{j}^{T}x))^{T},$$
}{\small
$$F_{i}(y) = (y[2i-1] - y[2i])^2.$$
}

\begin{table}[!htbp]

\caption{Experimental datasets for Reinforcement Learning}
\centering
\begin{tabular}{ c c c c } 
\hline
\hline
  &State number & Feature size & $\gamma$\\
\hline
Dataset-ref-1 & 2000 & 50 & $10^{-5}$ \\
Dataset-ref-2 & 2000 & 500 & $10^{-5}$ \\
\hline
\hline
\end{tabular}
\label{reinforce_data}
\end{table}

Following \citep{wang2016accelerating}, we generate a Markov decision problem (MDP) with totally 2000 states, and 10 actions for each state. The transition probability and the state features are randomly generated from the uniform distribution in the range of [0, 1]. In particular, we normalize the sum of transition probability from one state to 1. The dimension of state features is set to 50 for Dataset-ref-1 and 500 for Dataset-ref-2. The details of Dataset-ref-1 and Dataset-ref-2 can be found in Tabel~\ref{reinforce_data}. We add an $L_2$-regularization term $\frac{\gamma}{2} \| x \|^2$ to the loss function to make the strong convexity assumption hold and $\gamma$ is set to $10^{-5}$ to ensure the perturbation of the loss is small enough. AsyVRSC-Shared is run on this task with the number of threads varying from 1 to 16. We implement the following algorithms to compare with AsyVRSC-Shared:
\begin{itemize}
  \item SynSCGD-Shared: There are two phases in SCGD~\citep{wang2017stochastic}. In the first phase, all threads synchronously compute $y_t$. In the second phase, all threads synchronously compute the stochastic gradients and average them to conduct the update.
  \item SynVRSC-Shared: A synchronous parallel version of VRSC~\citep{lian2017finite}. In each inner iteration, all threads synchronously calculate the stochastic gradients and average them to conduct the update.
\end{itemize}

The results for AsyVRSC-Shared on this task are shown in Figure~\ref{figure_ref}. The objective value gap is defined as $|f(x)-f(x^*)|$. We draw the curves of objective value gap against time and iteration for AsynVRSC-Shared and compare AsyVRSC-Shared with SynSCGD-Shared and SynVRSC-Shared on the two datasets with 16 threads. From these results, we have the following observations: (i) AsynVRSC-Shared has linear convergence rate and more threads lead to less convergence time. (ii) AsyVRSC-Shared signigicantly outperforms SynSCGD-Shared and SynVRSC-Shared. (iii) AsyVRSC-Shared is more suitable for high dimensional problems than low dimensional problems.

\begin{figure}[!htb]
\centering

\subfigure[Dataset-port-1]{
\begin{minipage}[t]{0.3\linewidth}
\centering
\includegraphics[width=1.1\textwidth] {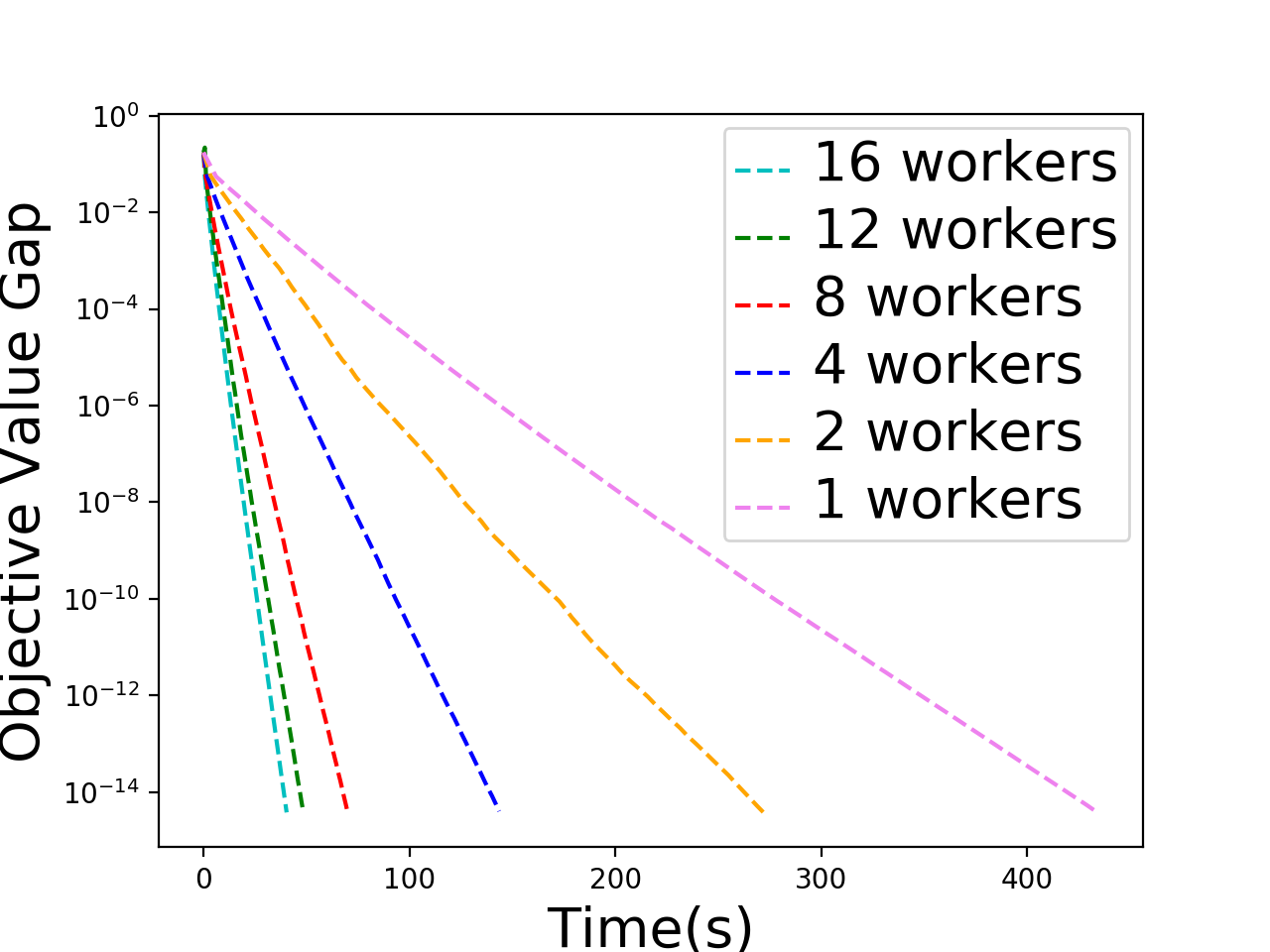}
\end{minipage}
}
\subfigure[Dataset-port-1]{
\begin{minipage}[t]{0.3\linewidth}
\centering
\includegraphics[width=1.1\textwidth] {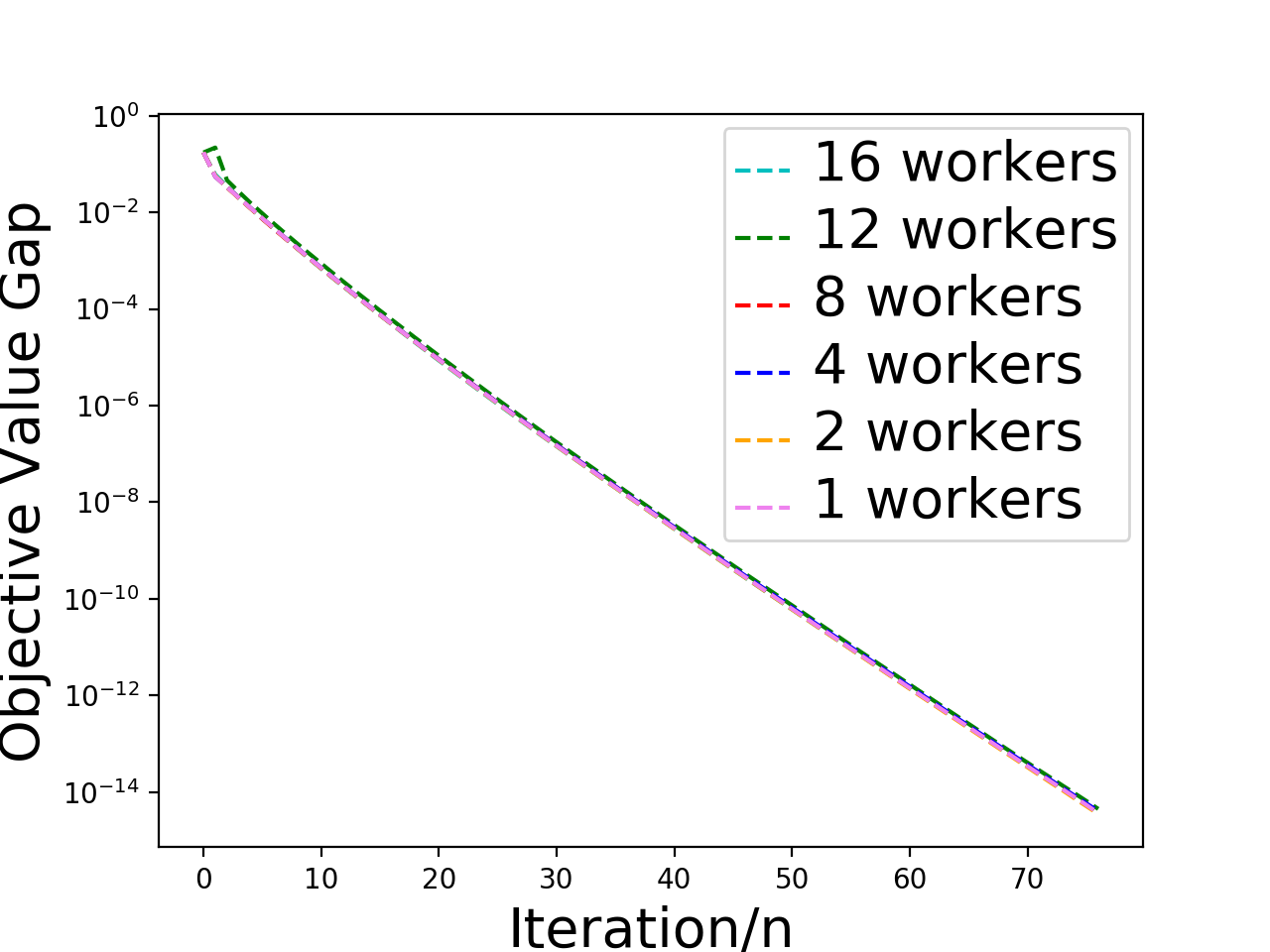}
\end{minipage}
}
\subfigure[Dataset-port-1]{
\begin{minipage}[t]{0.3\linewidth}
\centering
\includegraphics[width=1.1\textwidth] {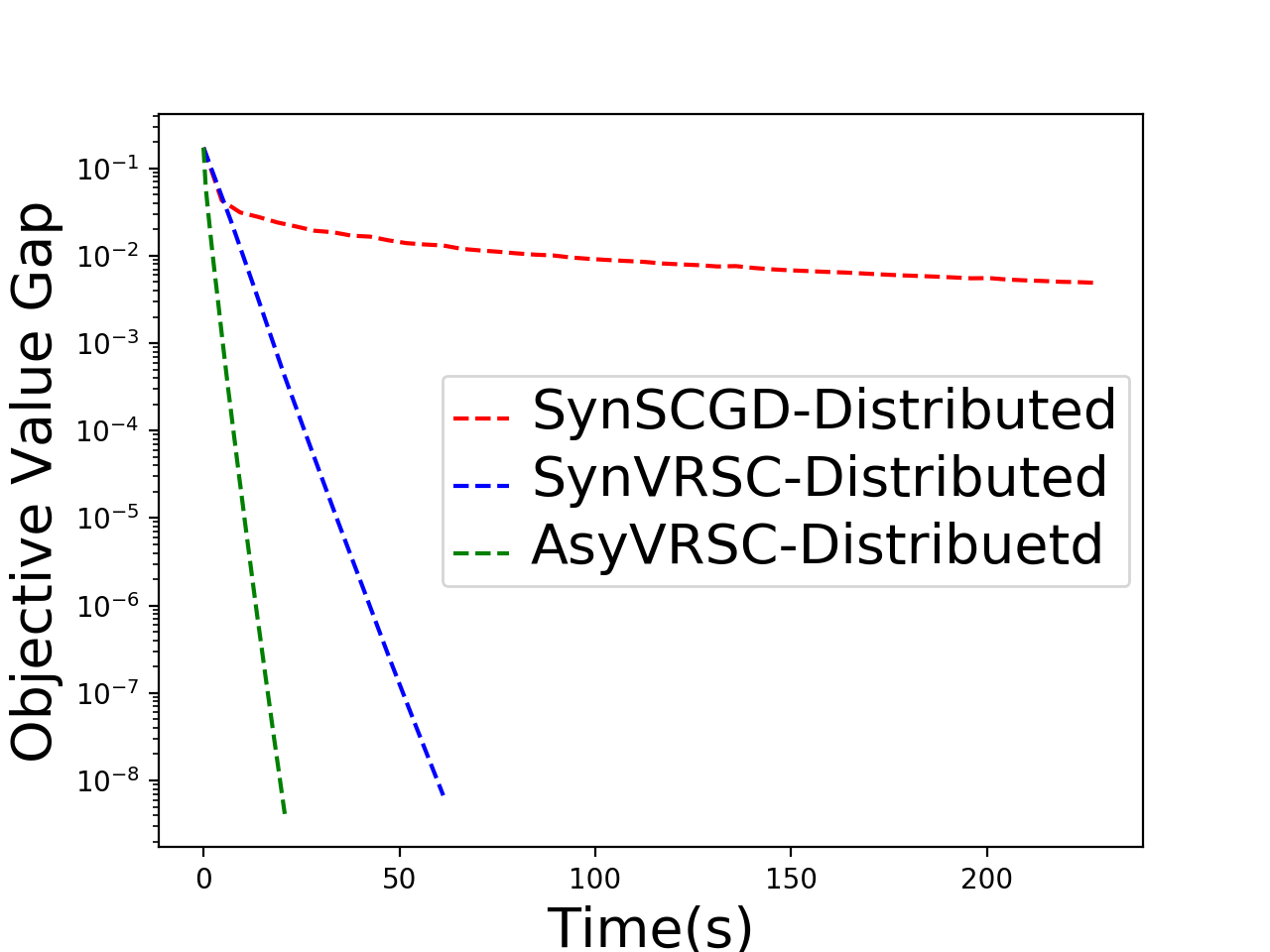}
\end{minipage}
}

\subfigure[Dataset-port-2]{
\begin{minipage}[t]{0.3\linewidth}
\centering
\includegraphics[width=1.1\textwidth] {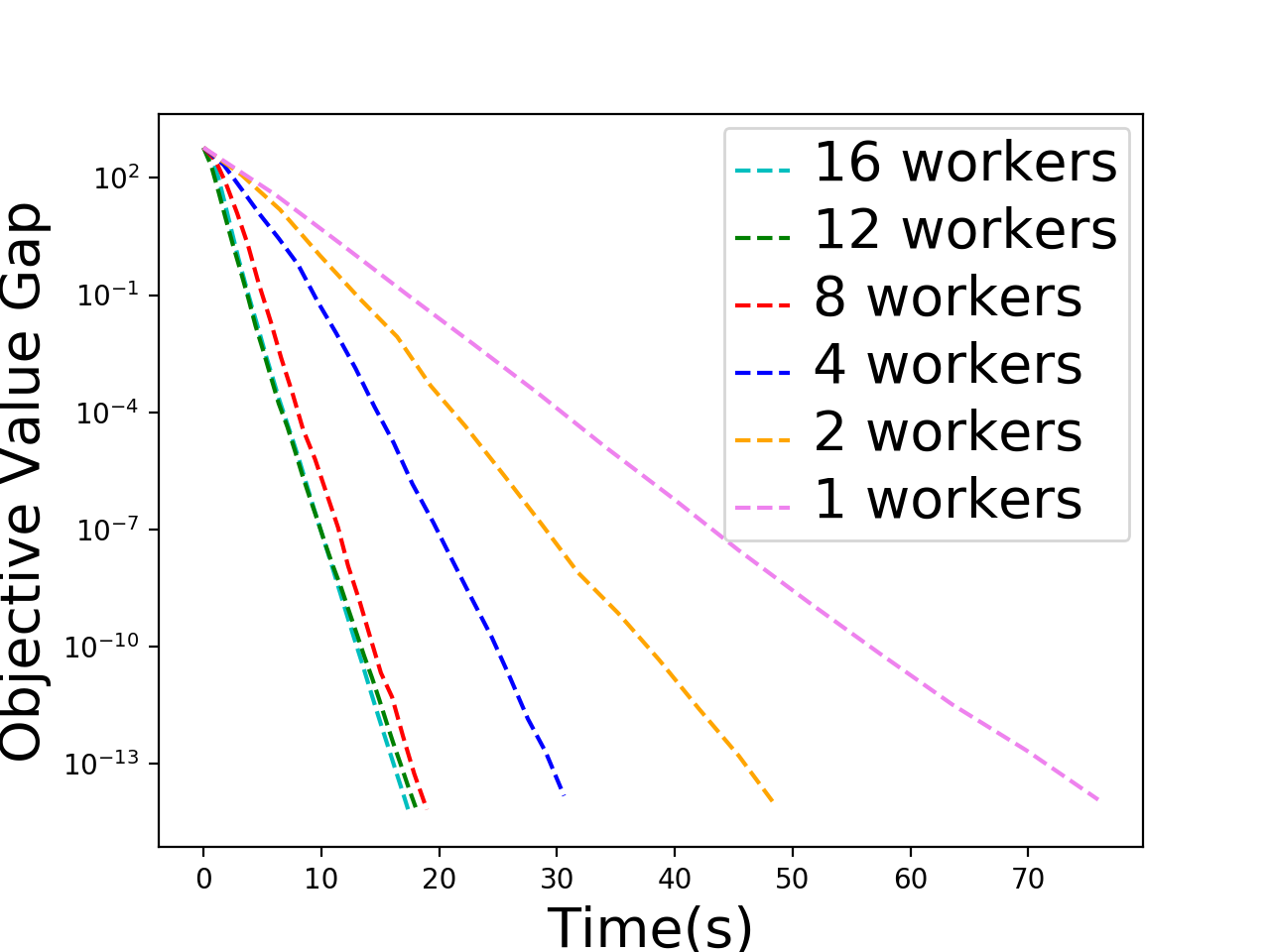}
\end{minipage}
}
\subfigure[Dataset-port-2]{
\begin{minipage}[t]{0.3\linewidth}
\centering
\includegraphics[width=1.1\textwidth] {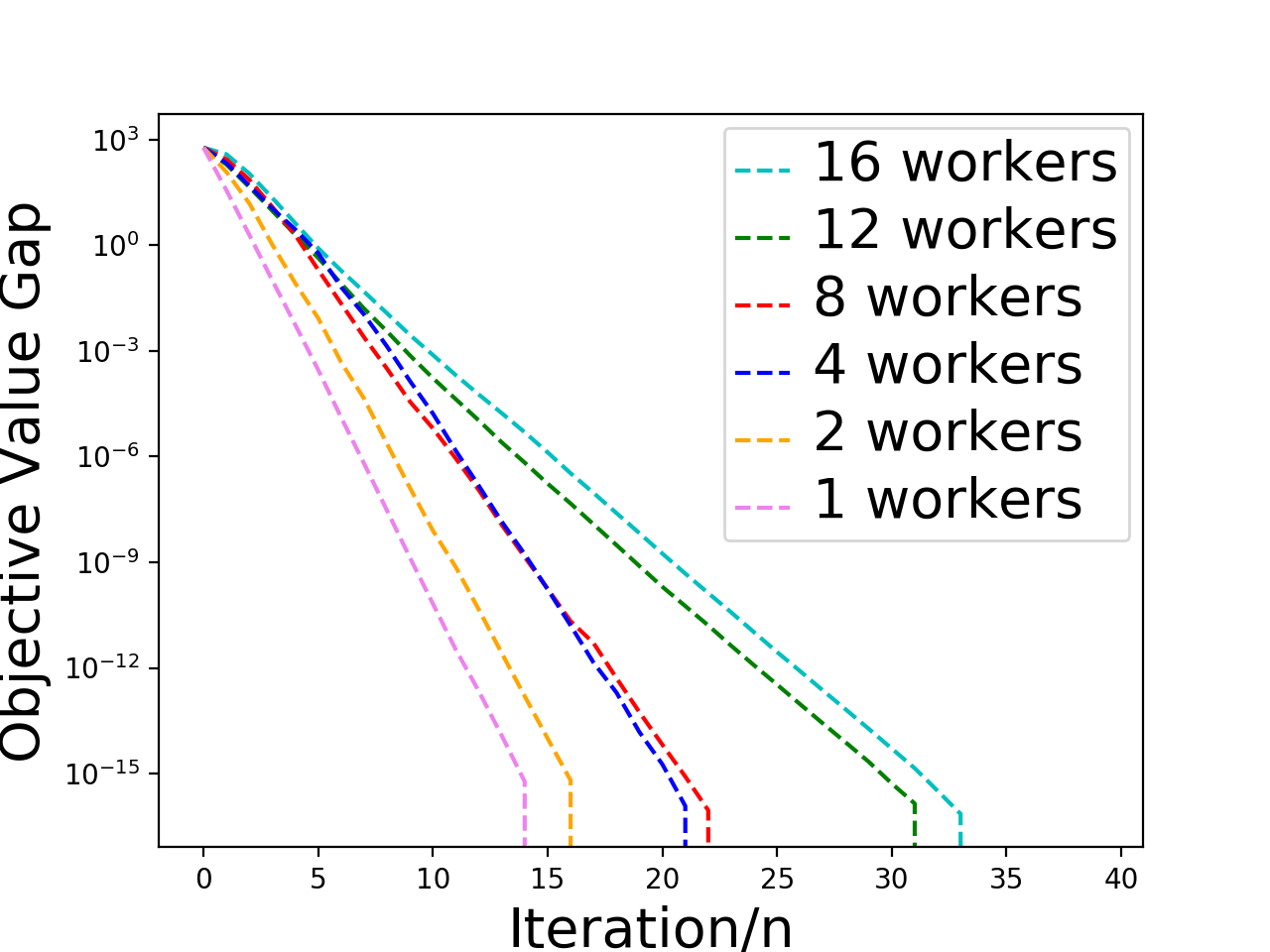}
\end{minipage}
}
\subfigure[Dataset-port-2]{
\begin{minipage}[t]{0.3\linewidth}
\centering
\includegraphics[width=1.1\textwidth] {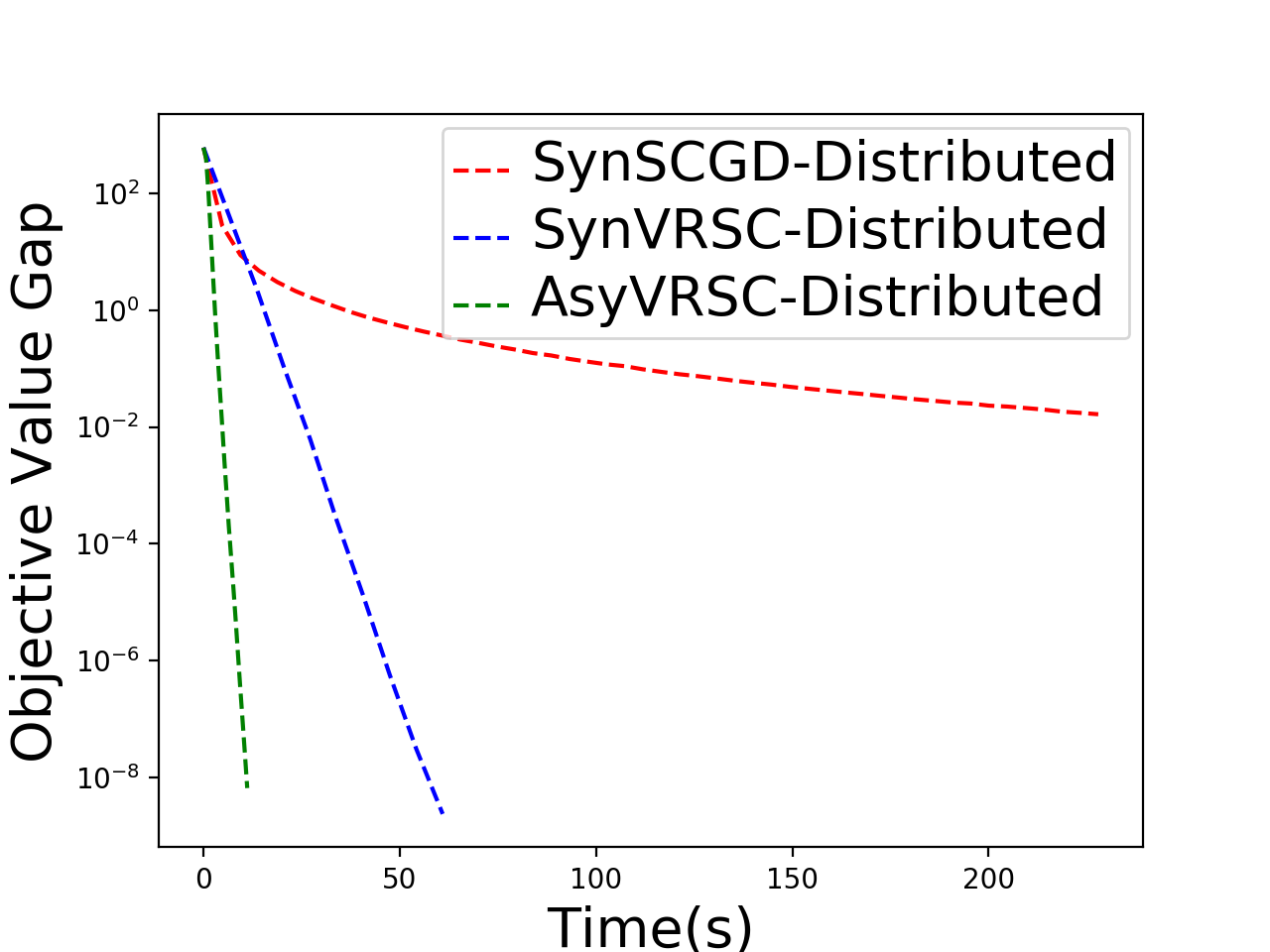}
\end{minipage}
}

\subfigure[Dataset-port-3]{
\begin{minipage}[t]{0.3\linewidth}
\centering
\includegraphics[width=1.1\textwidth] {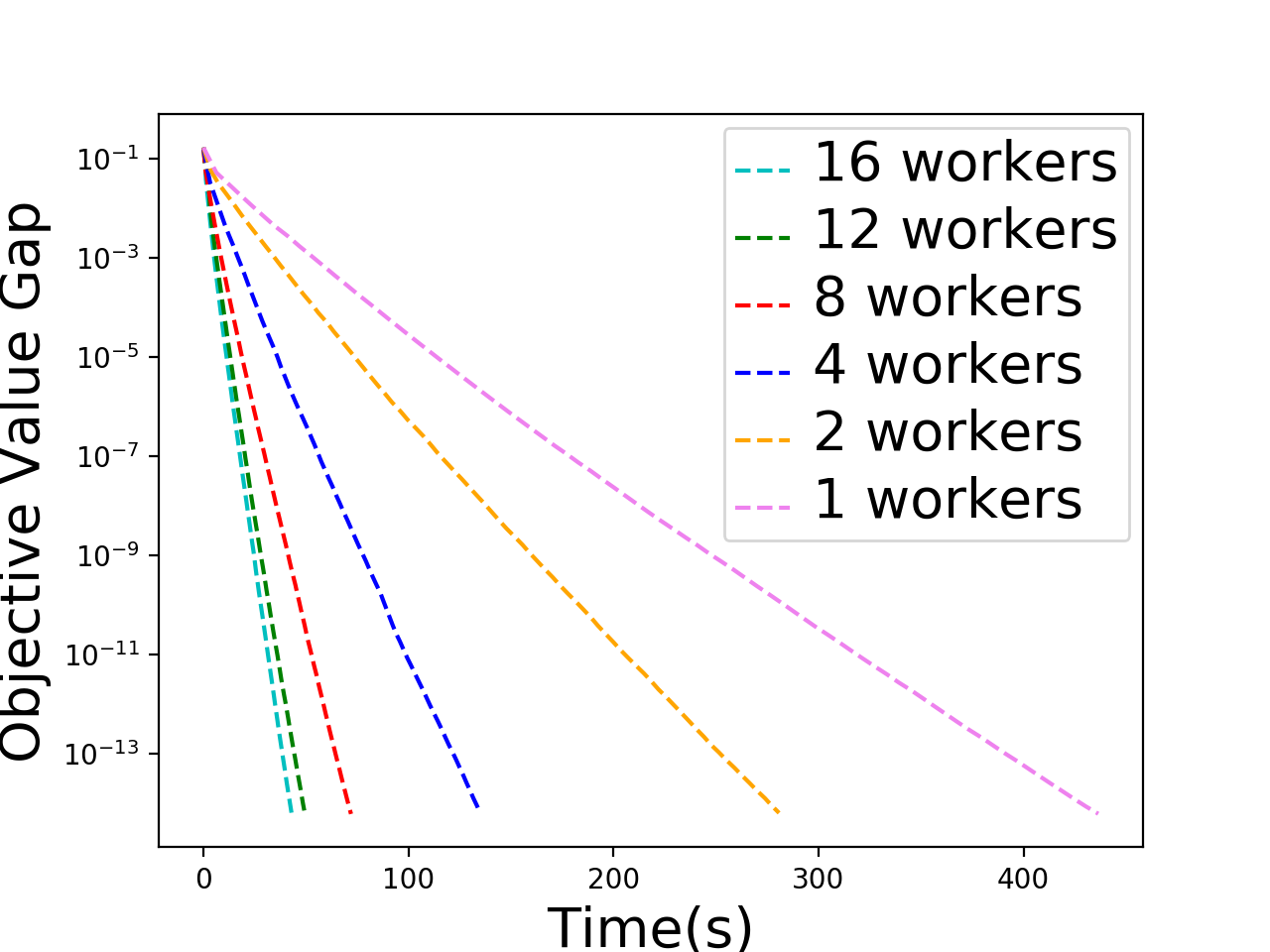}
\end{minipage}
}
\subfigure[Dataset-port-3]{
\begin{minipage}[t]{0.3\linewidth}
\centering
\includegraphics[width=1.1\textwidth] {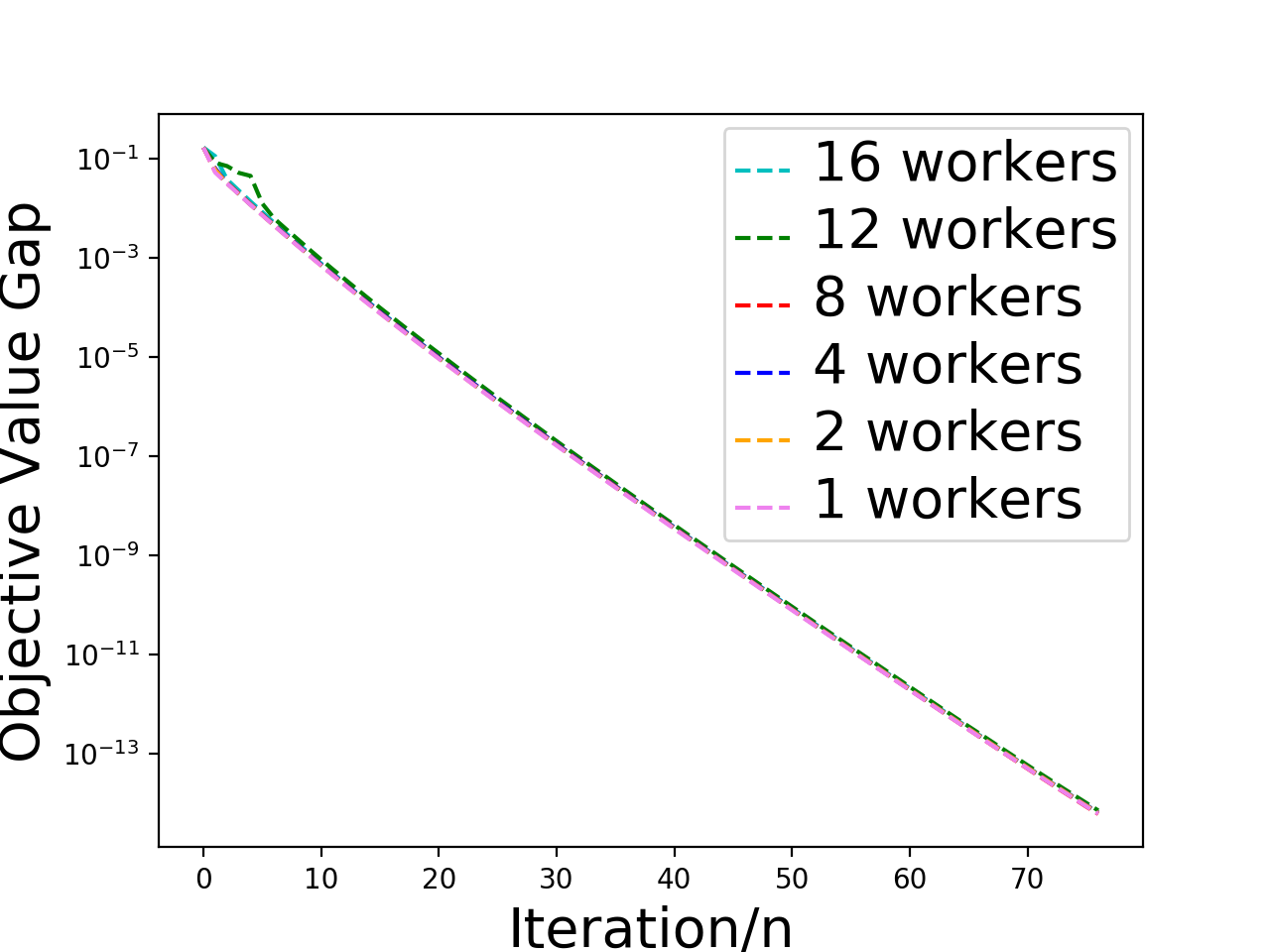}
\end{minipage}
}
\subfigure[Dataset-port-3]{
\begin{minipage}[t]{0.3\linewidth}
\centering
\includegraphics[width=1.1\textwidth] {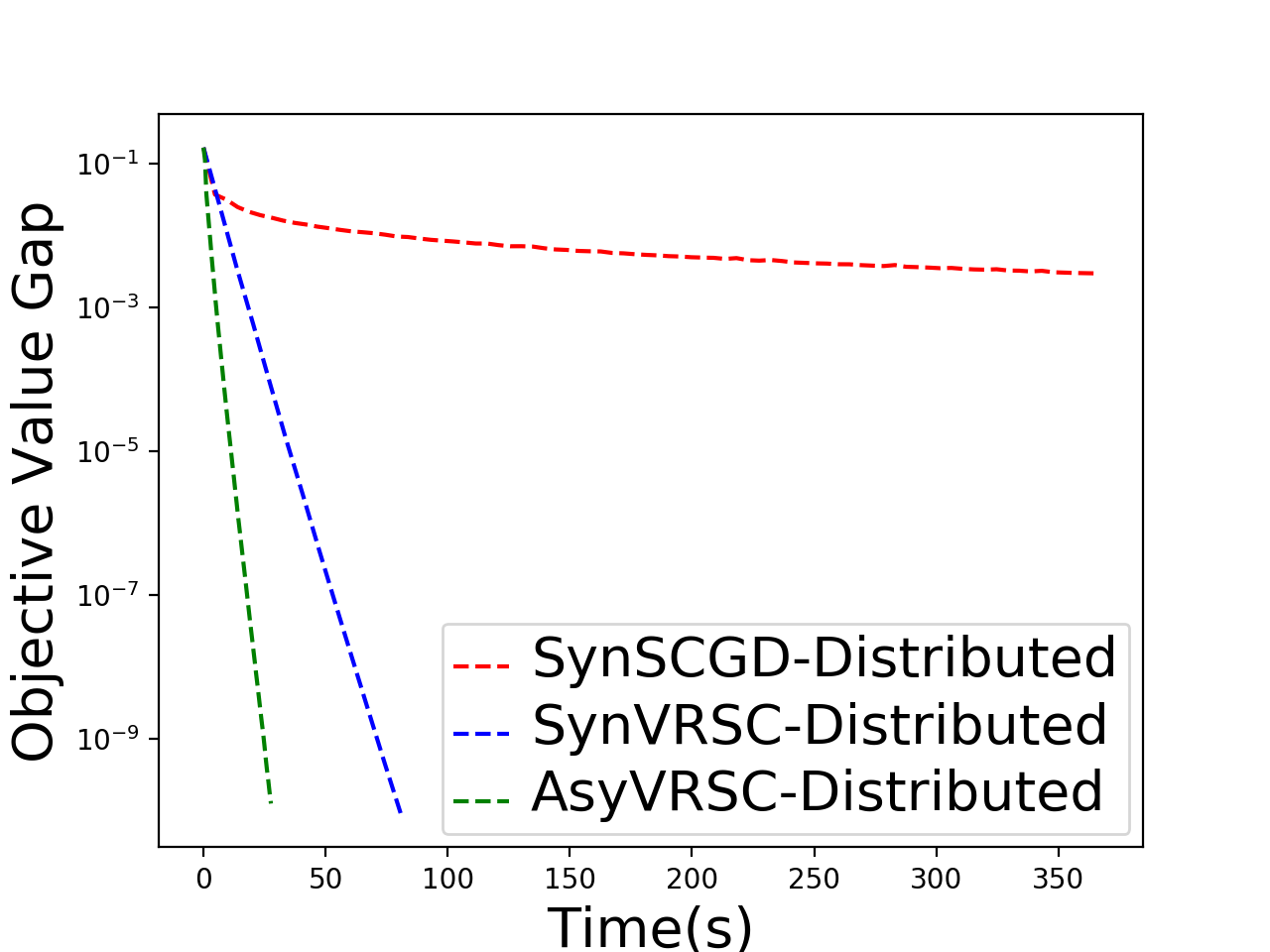}
\end{minipage}
}

\subfigure[Time Speedup]{
\begin{minipage}[t]{0.3\linewidth}
\centering
\includegraphics[width=1.1\textwidth] {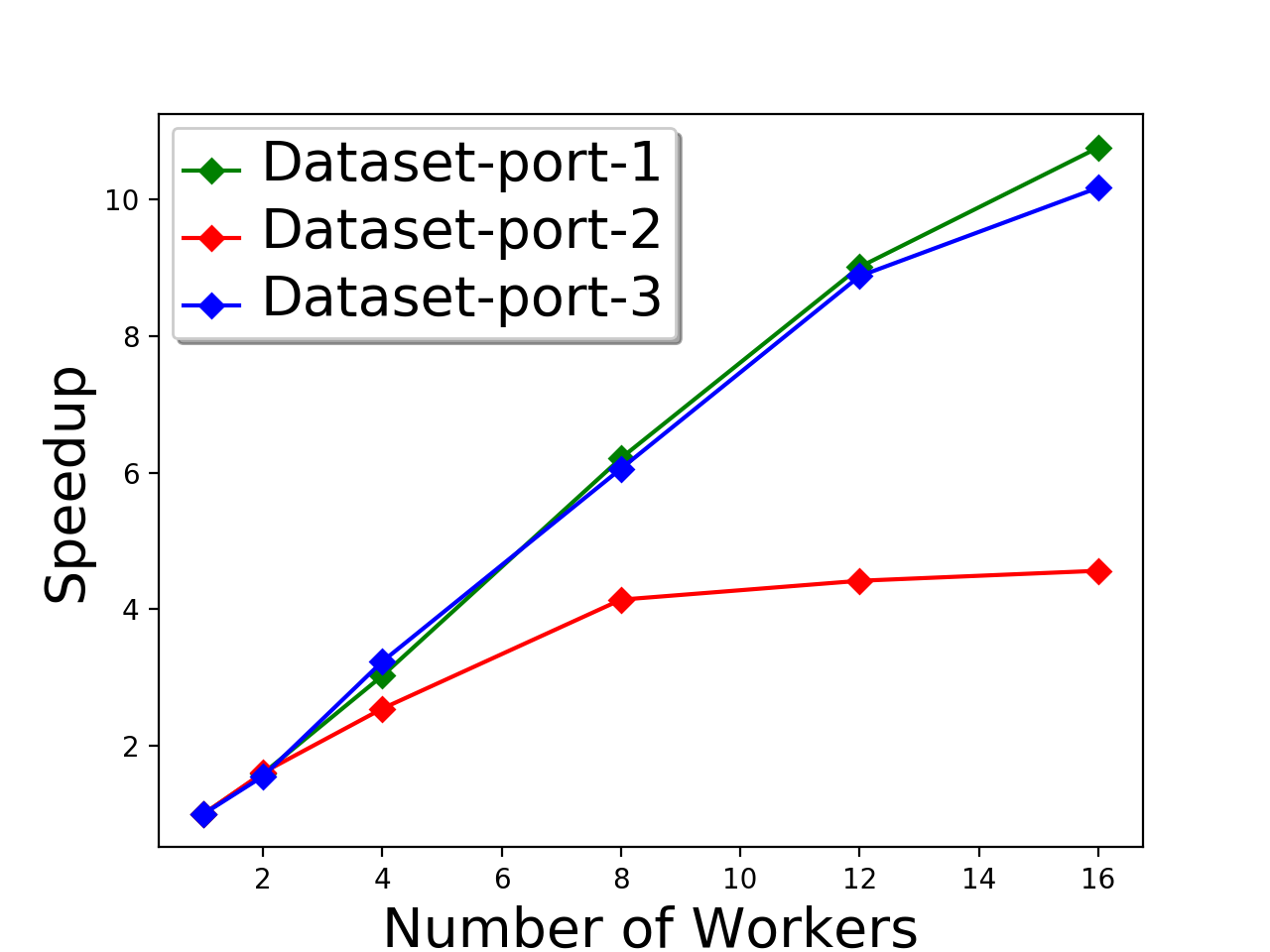}
\end{minipage}
}
\subfigure[Iteration Speedup]{
\begin{minipage}[t]{0.3\linewidth}
\centering
\includegraphics[width=1.1\textwidth] {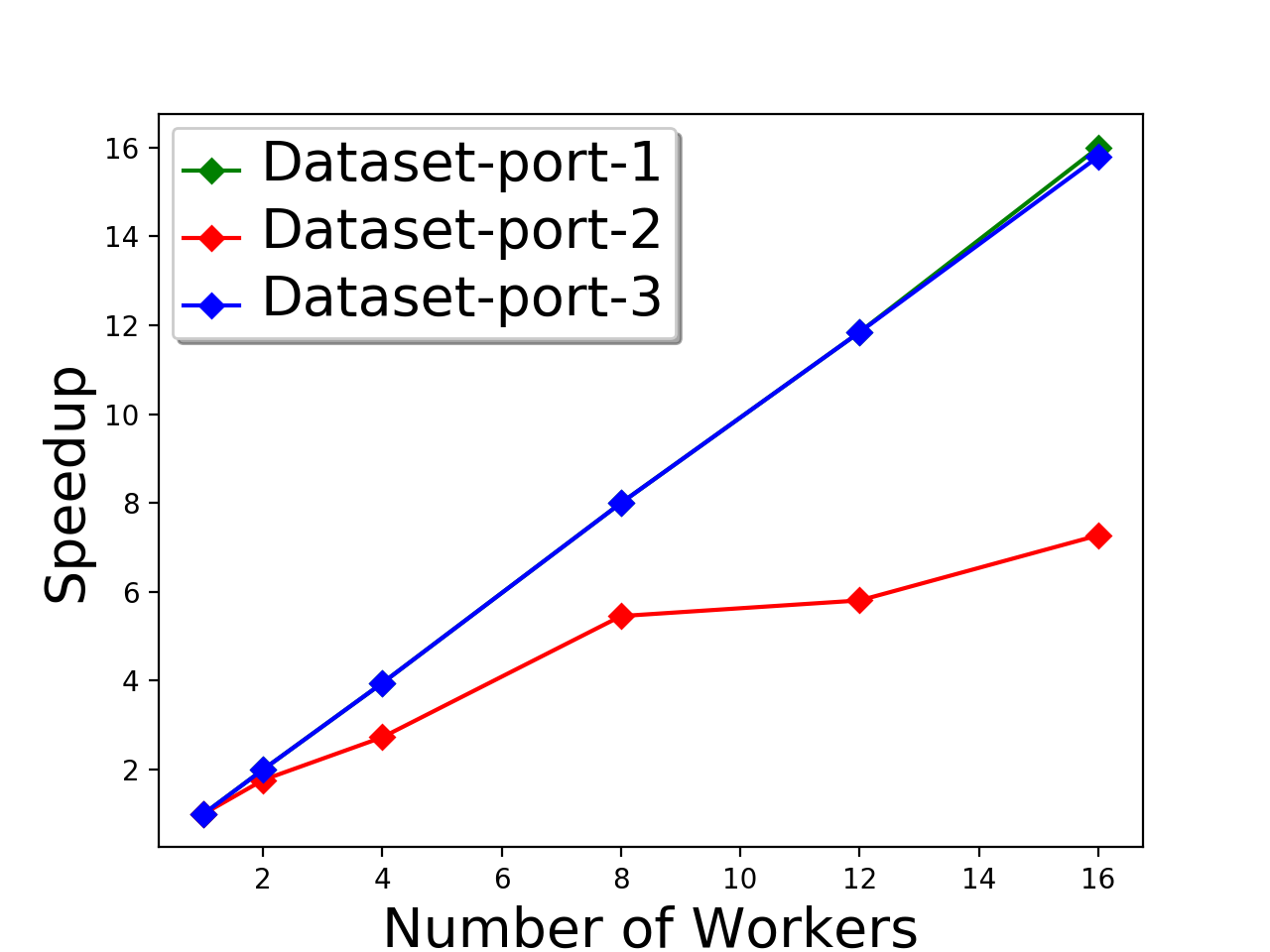}
\end{minipage}
}
 
\caption{Results for AsyVRSC-Distributed on the portfolio management task. Figure (a)-(c), (d)-(f) and (g)-(i) show the results for Dataset-port-1, Dataset-port-2 and Dataset-port-3, respectively. Figure (c), (f) and (i) show the comparision of AsyVRSC-Distribuetd with SynSCGD-Distribuetd and SynVRSC-Distribuetd with 16 workers on the three datasets. Figure (j) and (k) show the Time Speedup and Iteration Speedup of AsyVRSC-Distributed on the three datasets.}
\label{figure_port}
\end{figure}

\subsection{Portfolio Management}
For mean-variance optimization in portfolio management, suppose that there are $N$ assets we can invest and the reward vectors of the $N$ assets are denoted as $r_{t} \in \mathbb{R}^{N}(t=1,2,\cdots,n)$. The problem can then be formulated as:
{\small 
$$\min_{x \in \mathbb{R}^{N}} -\frac{1}{n} \sum_{i=1}^{n} \langle r_{i},x\rangle + \frac{1}{n} \sum_{i=1}^{n}\biggl(\langle r_{i},x \rangle - \frac{1}{n}\sum_{j=1}^{n} \langle r_{j},x \rangle \biggr)^2,$$
}where $x \in \mathbb{R}^{N}$ is the quantities invested to each portfolio. Same as~\citep{lian2017finite}, we use the following specifications for $G_{j}(x)$ and $F_{i}(y)$:
{\small
$$G_{j}(x) = (x^{T}, \langle r_{j},x \rangle)^{T},$$
$$F_{i}(y) = - y[N+1] + \biggl(\langle r_{i},y[1:N] \rangle - y[N+1]\biggr)^2.$$
}

\begin{table}[!htbp]

\caption{Experimental datasets for Portfolio Management}
\centering
\begin{tabular}{c c c c c c c } 
\hline
\hline
  &Data size & Feature size & $\lambda_{\max}$ & $\lambda_{\min}$ & $\gamma$ & sparsity \\
\hline
Dataset-port-1 & 10000 & 300 & 300 & 1 & 0 & 7.95\% \\
Dataset-port-2 & 10000 & 300 & 300 & 1 & 0 & 100\% \\
Dataset-port-3 & 10000 & 300 & 300 & 0 & $10^{-5}$ & 7.98\%\\
\hline
\hline
\end{tabular}
\label{port_data}
\end{table}

The reward vectors are generated in two steps: 
\begin{itemize}
  \item(1) Generate a vector in $R^{N}$ following the Gaussian distribution, where we define the maximum eigenvalue and the minimum eigenvalue of its covariance matrix as $\lambda_{\max}$ and $\lambda_{\min}$, respectively. Because the condition number of its covariance matrix $\frac{\lambda_{\max}}{\lambda_{\min}}$ is proportional to $\kappa$, we will use $\lambda_{\max}$ and $\lambda_{\min}$ to control the Lipschitz gradient and the strong convexity constants defined in Assumption~\ref{assumption:smoothness} and Assumption~\ref{assumption:strongly convex}. 
  \item(2) Sample rewards $r_t$ from the Gaussian distribution and set all elements to its absolute value to ensure the problem has an optimal solution. 
\end{itemize}
We generate three datasets with the prarameters shown in Table~\ref{port_data}. In the experiment, we add an $L_2$-regularization term $\frac{\gamma}{2} \| x \|^2$ to the loss function of Dataset-port-3 to make the strong convexity assumption hold and $\gamma$ is set to $10^{-5}$ to ensure the perturbation of the loss small enough. AsyVRSC-Distributed is run on this task with the number of workers varying from 1 to 16. We implement SynSCGD-Distributed and SynVRSC-Distributed similar as section~\ref{section_ref} to compare with AsyVRSC-Distributed.

The results are demonstrated in Figure~\ref{figure_port}. We draw the curves of objective value gap against time and iteration for AsyVRSC-Distributed and compare AsyVRSC-Distribuetd to SynSCGD-Distributed and SynVRSC-Distributed  
with 16 workers. For Dataset-port-1 and Dataset-port-3, $\nabla F_{i}(x)$, $\nabla G_{j}(x)$ and $f_{ij}(x)$ satisfy Assumption 5. From these results, we have the following observations:  
(i) AsyVRSC-Distribuetd has linear convergence rate and more workers lead to less convergence time.
(ii) AsynVRSC-Distribuetd can significantly outperform SynSCGD-Distribuetd and SynVRSC-Distributed.
(iii) AsyVRSC-Distributed achieves linear speedup when the sparsity assumption holds and the perfermance of AsyVRSC-Distribuetd in sparse problems is better than that in the dense problems.

\section{Conclusion}
In this paper, we study the asynchronous parallelization of stochastic composition optimization with variance reduction. We propose AsyVRSC-Shared and AsyVRSC-Distributed for shared-memory architecture and master-worker architecture, respectively. We prove that both proposed algorithms can achieve linear convergence rate for strongly convex loss functions. When the number of workers grows, both algorithms have linear speedup under certain conditions. Experimental results demonstrate the efficiency of the proposed algorithms.

\allowdisplaybreaks
\section*{Appendix A: Some Basic Lemmas}

\begin{lemma}
For any $\alpha > 0$ and $x, y \in \mathbb{R}^{d}$, we have
\begin{equation}\label{young's2}
-2 \alpha \|x\|^2 - \frac{1}{2 \alpha}\|y\|^2 \leq \langle x, y \rangle \leq 2 \alpha \|x\|^2 + \frac{1}{2 \alpha}\|y\|^2.
\end{equation} 
\end{lemma}

\begin{lemma}
For any $x, y \in \mathbb{R}^{d}$, we have
\begin{equation}\label{cauthy}
\langle x, y \rangle  \leq \| x \| \| y \| \leq \frac{1}{2} (\| x \|^2 + \| y \|^2).
\end{equation}
\end{lemma}

\begin{lemma}
For any variables $\beta_{1}, \cdots, \beta_{t} \in \mathbb{R}^d$, we have
\begin{equation}\label{young's}
\| \beta_{1} + \beta_2+\cdots+ \beta_{t} \|^2 \leq t (\| \beta_1\|^2 + \| \beta_2 \|^2 + \cdots + \| \beta_t\|^2), \forall t \in \mathbb{N}_{+}.
\end{equation}
\end{lemma}

\begin{lemma}
For any random vector $\zeta \in \mathbb{R}^{d}$, it holds that 
\begin{equation} \label{expectation}
\mathbb{E} \| \zeta - \mathbb{E} \zeta \|^2 = \mathbb{E} \| \zeta \|^2 - \| \mathbb{E} \zeta \|^2.
\end{equation}
\end{lemma}

\section*{Appendix B: Convergence Analysis for Section 4.1}

\textbf{Proof of Lemma~\ref{lemma:1}.}

\noindent We define an unbiased estimation of $\nabla f(x_{t}^{s})$ as
$$\nabla \widetilde{f}(x_{t}^{s}) = (\nabla G(x_{t}^{s}))^{T} \nabla F_{i_t}(G(x_{t}^{s})) - \nabla f_{i_t}(\widetilde{x}^{s}) + \nabla f(\widetilde{x}^{s}).$$ 
First, we bound $\| \nabla \widetilde{f}(x_{t}^{s}) - \nabla \widehat{f}(x_{t}^{s}) \|^2$ by
\begin{eqnarray} \label{lemma_1_1}
&&\mathbb{E} \|\nabla \widetilde{f}(x_{t}^{s}) - \nabla \widehat{f}(x_{t}^{s})\|^{2}  
\nonumber\\
& = & \mathbb{E} \| (\nabla G(x_{t}^{s}))^{T} \nabla F_{i_t}(G(x_{t}^{s})) - (\nabla \widehat{G}_t^{s})^{T} \nabla F_{i_t}(\widehat{G}_{t}^{s})\|^{2}
\nonumber\\
&=& \mathbb{E} \| (\nabla G(x_{t}^{s}))^{T} \nabla F_{i_t}(G(x_{t}^{s}))  - (\nabla G(x_{t}^{s}))^T \nabla F_{i_{t}}( \widehat{G}_{t}^{s}) + (\nabla G(x_{t}^{s}))^T \nabla F_{i_{t}}( \widehat{G}_{t}^{s}) 
\nonumber\\
&&- (\nabla \widehat{G}_t^{s})^{T} \nabla F_{i_t}(\widehat{G}_{t}^{s})\|^{2}
\nonumber\\
& \overset{(\ref{young's})}{\leq}& 2\mathbb{E}  \|(\nabla G(x_{t}^{s}))^{T} \nabla F_{i_t}(G(x_{t}^{s})) - (\nabla G(x_{t}^{s}))^T \nabla F_{i_{t}}( \widehat{G}_{t}^{s} ) \|^2  
\nonumber\\
&& + 2\mathbb{E} \|(\nabla G(x_{t}^{s}))^T \nabla F_{i_{t}}( \widehat{G}_{t}^{s}) - (\nabla \widehat{G}_t^{s})^{T} \nabla F_{i_t}(\widehat{G}_{t}^{s})\|^{2} 
\nonumber\\
& \overset{(\ref{cauthy})}{\leq}& 2\mathbb{E} \| \nabla G(x_{t}^{s}) \|^2 \| \nabla F_{i_t}(G(x_{t}^{s})) -  \nabla F_{i_{t}}( \widehat{G}_{t}^{s}) \|^2 + 2\mathbb{E} \|\nabla G(x_{t}^{s}) 
\nonumber\\
&&- \nabla \widehat{G}_{t}^{s}\|^2 \| \nabla F_{i_t}(\widehat{G}_{t}^{s}) \|^2
\nonumber\\
&\leq& 2B_{G}^2 \mathbb{E} \|\nabla F_{i_t}(G(x_{t}^{s})) - \nabla F_{i_t}(\widehat{G}_{t}^{s})\|^2 + 2 B_{F}^2 \mathbb{E} \|\nabla G(x_{t}^{s}) - \nabla \widehat{G}_{t}^{s}\|^2
\nonumber\\
&\leq& 2B_{G}^2 L_{F}^2 \mathbb{E} \|G(x_{t}^{s}) - G(\widetilde{x}^{s}) - \frac{1}{a} \sum_{j=1}^{a}(G_{A_{t}[j]}(x_{t}^{s}) - G_{A_{t}[j]}(\widetilde{x}^{s}))\|^2
\nonumber\\
&&+ 2B_{F}^2 \mathbb{E} \|\nabla G(x_{t}^{s}) - \nabla G(\widetilde{x}^s) - \frac{1}{b} \sum_{j=1}^{b} ( \nabla G_{B_{t}[j]}(x_{t}^{s}) - \nabla G_{B_{t}[j]}( \widetilde{x}^{s}))\|^2
\nonumber\\
&=& \frac{2B_{G}^2 L_{F}^2}{a^2} \mathbb{E} \| \sum_{j=1}^{a} (G(x_{t}^{s}) - G(\widetilde{x}^{s}) - (G_{A_{t}[j]}(x_{t}^{s}) - G_{A_{t}[j]}(\widetilde{x}^{s})) )\|^2
\nonumber\\
&& + \frac{2B_{F}^2}{b^2} \mathbb{E} \|\sum_{j=1}^{b}(\nabla G(x_{t}^{s}) - \nabla G(\widetilde{x}^s) - (\nabla G_{B_{t}[j]}(x_{t}^{s}) - \nabla G_{B_{t}[j]}( \widetilde{x}^{s})) )\|^2
\nonumber\\
&= & \frac{2B_{G}^2 L_{F}^2}{a^2} \sum_{j=1}^{a} \mathbb{E} \|G(x_{t}^{s}) - G(\widetilde{x}^{s}) - (G_{A_{t}[j]}(x_{t}^{s}) - G_{A_{t}[j]}(\widetilde{x}^{s}))\|^2
\nonumber\\
& & + \frac{2B_{F}^2}{b^2} \sum_{j=1}^{b} \mathbb{E} \|\nabla G(x_{t}^{s}) - \nabla G(\widetilde{x}^s) - (\nabla G_{B_{t}[j]}(x_{t}^{s}) - \nabla G_{B_{t}[j]}( \widetilde{x}^{s}))\|^2
\nonumber\\
& \overset{(\ref{young's})}{\leq}& \frac{8B_{G}^2 L_{F}^2}{a^2} \mathbb{E} \sum_{j=1}^{a}(\|G(x_{t}^{s}) - G(x^{*})\|^2 + \|G(\widetilde{x}^{s}) - G(x^{*})\|^2 + \|G_{A_{t}[j]}(x_{t}^{s}) - G_{A_{t}[j]}(x^{*})\|^2 
\nonumber\\
& & + \|G_{A_{t}[j]}(\widetilde{x}^{s}) - G_{A_{t}[j]}(x^{*})\|^2) + \frac{8B_{F}^2}{b^2} \mathbb{E} \sum_{j=1}^{B} (\|\nabla G(x_{t}^{s}) - \nabla G(x^{*})\|^2 
\nonumber\\
& & + \|\nabla G(\widetilde{x}^{s}) - \nabla G(x^{*})\|^2 + \|\nabla G_{B_{t}[j]}(x_{t}^{s}) - \nabla G_{B_{t}[j]}(x^{*})\|^2 
\nonumber\\
&&+ \|\nabla G_{B_{t}[j]}(\widetilde{x}^{s}) - \nabla G_{B_{t}[j]}(x^{*})\|^2),
\end{eqnarray}
where the third inequality and the fourth inequality follow from Assumption \ref{assumption:bounded gradient} and Assumption \ref{assumption:smoothness}, respectively. The last equality comes from the fact that the indices in $A_{t}$ and $B_{t}$ are independent. Specifically,
$\mathbb{E}_{i \neq j} \langle G(x_{t}^{s}) - G^{s} - (G_{A_{t}[i]}(x_{t}^{s}) - G_{A_{t}[i]}(\widetilde{x}^{s})), G(x_{t}^{s}) - G^{s} - (G_{A_{t}[j]}(x_{t}^{s}) - G_{A_{t}[j]}(\widetilde{x}^{s})) \rangle = 0.$
Combining Assumption \ref{assumption:bounded gradient} and Assumption \ref{assumption:smoothness} with (\ref{lemma_1_1}), $\| \nabla \widetilde{f}(x_{t}^{s}) - \nabla \widehat{f}(x_{t}^{s}) \|^2$ can be finally bounded by
\begin{eqnarray}\label{lemma_1_2}
&&\mathbb{E} \|\nabla \widetilde{f}(x_{t}^{s}) - \nabla \widehat{f}(x_{t}^{s})\|^{2}
\nonumber\\
&\leq& \frac{8B_{G}^2 L_{F}^2}{a^2} \mathbb{E} \sum_{j=1}^{a} 2B_{G}^2 (\|x_{t}^{s} - x^{*}\|^2 + \|\widetilde{x}^{s} - x^{*}\|^2) 
\nonumber\\
&&+ \frac{8B_{F}^2}{b^2} \mathbb{E} \sum_{j=1}^{b} 2L_{G}^2 (\|x_{t}^{s} - x^{*}\|^2 + \|\widetilde{x}^{s} - x^{*}\|^2)
\nonumber\\
&=& 16(\frac{B_{G}^4 L_{F}^2}{a} + \frac{B_{F}^2 L_{G}^2}{b}) \mathbb{E} (\|x_{t}^{s} - x^{*}\|^2 + \|\widetilde{x}^{s} - x^{*}\|^{2})
\nonumber\\
&\leq& \frac{32}{\mu_{f}} (\frac{B_{G}^4 L_{F}^2}{a} + \frac{B_{F}^2 L_{G}^2}{b}) \mathbb{E}(f(x_{t}^{s}) - f(x^{*}) + f(\widetilde{x}^{s}) - f(x^{*})),
\end{eqnarray}
where the last inequality follows from that $f(x)$ is $\mu_{f}$-strongly convex. Next, we bound $\| \nabla \widetilde{f}(x_{t}^{s}) \|$ by
\begin{eqnarray}\label{lemma_1_3}
&&\mathbb{E} \| \nabla \widetilde{f}(x_{t}^{s}) \|^{2} 
\nonumber\\
& = & \mathbb{E}\|(\nabla G(x_{t}^{s}))^{T} \nabla F_{i_{t}}(G(x_{t}^{s})) - (\nabla G(\widetilde{x}))^{T} \nabla F_{i_{t}}(G^{s}) + \nabla f(\widetilde{x}^{s})\|^2
\nonumber\\
& \overset{(\ref{young's})}{\leq}& 2 \mathbb{E} \|(\nabla G(x_{t}^{s}))^{T} \nabla F_{i_{t}}(G(x_{t}^{s})) - (\nabla G(x^{*}))^{T} \nabla F_{i_{t}}(G(x^{*}))\|^{2}
\nonumber\\
& & + 2 \mathbb{E} \|(\nabla G(\widetilde{x}) \nabla F_{i_{t}}(G^{s})) - (\nabla G(x^{*}))^{T} \nabla F_{i_{t}}(G(x^{*})) - (\nabla f(\widetilde{x}^{s}) - \nabla f(x^{*}))\|^2
\nonumber\\
&\overset{(\ref{expectation})}{=}& 2 \mathbb{E} \|(\nabla G(x_{t}^{s}))^{T} \nabla F_{i_{t}}(G(x_{t}^{s})) - (\nabla G(x^{*}))^{T} \nabla F_{i_{t}}(G(x^{*}))\|^{2}
\nonumber\\
& & + 2 \mathbb{E} \|(\nabla G(\widetilde{x}) \nabla F_{i_{t}}(G^{s})) - (\nabla G(x^{*}))^{T} \nabla F_{i_{t}}(G(x^{*}))\|^2 - 2\|\nabla f(\widetilde{x}^{s}) - \nabla f(x^{*})\|^2
\nonumber\\
&\leq& 2 \mathbb{E} \|(\nabla G(x_{t}^{s}))^{T} \nabla F_{i_{t}}(G(x_{t}^{s})) - (\nabla G(x^{*}))^{T} \nabla F_{i_{t}}(G(x^{*}))\|^{2}
\nonumber\\
& & + 2 \mathbb{E} \|(\nabla G(\widetilde{x}) \nabla F_{i_{t}}(G^{s})) - (\nabla G(x^{*}))^{T} \nabla F_{i_{t}}(G(x^{*}))\|^2
\nonumber\\
&\leq& 4 L_{f} ( (f(x_{t}^{s}) - f(x^{*}) + f(\widetilde{x}^{*}) - f(x^{*})),
\end{eqnarray}
where the last inequality comes from the smoothness assumption: $\| \nabla f_{i}(x) - \nabla f_{i}(x^*) \| \leq 2 L_{f} [f_{i}(x) - f_{i}(x^*) - \langle \nabla f_{i}(x), x- x^* \rangle ]$ (Theorem 2.1.5 in \citep{nesterov2013introductory}). 
Combining (\ref{lemma_1_2}) and (\ref{lemma_1_3}), we can finally get
\begin{eqnarray} \label{lamma_1_4}
&&\mathbb{E} \|\nabla \widehat{f}(x_{t}^{s})\|^2 
\nonumber\\
& \overset{(\ref{young's})}{\leq}& 2\mathbb{E} \|\nabla \widehat{f}(x_{t}^{s}) - \nabla \widetilde{f}(x_{t}^{s}) \|^{2} + 2 \mathbb{E} \| \nabla \widetilde{f}(x_{t}^{s}) \|^{2}
\nonumber\\
&\leq& (\frac{64}{\mu_{f}}(\frac{B_{G}^4 L_{F}^2}{a} + \frac{B_{F}^2 L_{G}^2}{b}) + 8 L_{f}) (f(x_{t}^{s}) - f(x^{*}) + f(\widetilde{x}^{s}) - f(x^{*})).
\end{eqnarray}

\noindent \textbf{Proof of Lemma~\ref{lemma:2}.}

\noindent Based on the definition of $\widehat{G}_{t - \tau_{t}^{s}}^{s}$, we have
\begin{eqnarray}\label{shared_T4_1}
&&\|\widehat{G}_{t - \tau_{t}^{s}}^{s} - G(x_{t - \tau_{t}^{s}}^{s})\|^2 
\nonumber\\
&=& \mathbb{E} \|G(\widetilde{x}^{s}) - \frac{1}{a} \sum_{j=1}^{a}( G_{A_{t}[j]}(\widetilde{x}^{s}) - G_{A_{t}[j]}(x_{t - \tau_{t}^{s}}^{s}) ) - G(x_{t - \tau_{t}^{s}}^{s})\|^2
\nonumber\\
&=& \frac{1}{a^2} \mathbb{E} \|\sum_{j=1}^{a}(G_{A_{t}[j]}(\widetilde{x}^{s}) - G_{A_{t}[j]}(x_{t - \tau_{t}^{s}}^{s}) - (G(\widetilde{x}^{s}) - G(x_{t - \tau_{t}^{s}}^{s})))\|^2
\nonumber\\
&=& \frac{1}{a^2} \sum_{j=1}^{a} \mathbb{E} \|G_{A_{t}[j]}(\widetilde{x}^{s}) - G_{A_{t}[j]}(x_{t - \tau_{t}^{s}}^{s}) - (G(\widetilde{x}^{s}) - G(x_{t - \tau_{t}^{s}}^{s}))\|^2
\nonumber\\
&\overset{(\ref{expectation})}{=}& \frac{1}{a^2} \sum_{j=1}^{a} (\mathbb{E} \|G_{A_{t}[j]}(\widetilde{x}^{s}) - G_{A_{t}[j]}(x_{t - \tau_{t}^{s}}^{s})\|^2 - \mathbb{E} \|G(\widetilde{x}^{s}) - G(x_{t - \tau_{t}^{s}}^{s})\|^2)
\nonumber\\
&\leq& \frac{1}{a^2} \sum_{j=1}^{a} \mathbb{E} \|G_{A_{t}[j]}(\widetilde{x}^{s}) - G_{A_{t}[j]}(x_{t - \tau_{t}^{s}}^{s})\|^2 
\nonumber\\
&\leq& \frac{B_{G}^2}{a} \mathbb{E} \|\widetilde{x}^{s} - x_{t - \tau_{t}^{s}}^{s}\|^2
\nonumber\\
&=& \frac{B_{G}^2}{a} \mathbb{E} \|\widetilde{x}^{s} - x^{*} -x_{t}^{s} + x^{*} + x_{t}^{s}- x_{t - \tau_{t}^{s}}^{s}\|^2
\nonumber\\
&\overset{(\ref{young's})}{\leq}& \frac{3B_{G}^2}{a} \mathbb{E} (\|\widetilde{x}^{s} - x^{*}\|^2 + \|x_{t}^{s} - x^{*}\|^2 + \|x_{t}^{s}- x_{t - \tau_{t}^{s}}^{s}\|^2) \label{shared_T4_2},
\end{eqnarray}
where the third equality comes from the fact that the indices in $A_t$ are independent. The second inequality comes from Assumption~\ref{assumption:bounded gradient}. \\

\noindent \textbf{Proof for Lemma~\ref{lemma:3}.}

\noindent Based on the definition of $\nabla \widehat{G}_{t - \tau_{t}^{s}}^{s}$, we have
\begin{eqnarray}\label{shared_T6_1}
&&\mathbb{E}\|\nabla \widehat{G}_{t - \tau_{t}^{s}}^{s}\|^2
\nonumber\\
&=& \mathbb{E}\|\nabla G(\widetilde{x}^s) - \frac{1}{b} \sum_{j=1}^{b}(\nabla G_{B_t[j]}(\widetilde{x}^s) - \nabla G_{B_{t}[j]}(x_{t - \tau_{t}^{s}}^{s}) ) \|^2
\nonumber\\
&\overset{(\ref{young's})}{\leq}& 2 \mathbb{E}\|\nabla G(\widetilde{x}^s) - \frac{1}{b} \sum_{j=1}^{b}(\nabla G_{B_t[j]}(\widetilde{x}^s) - \nabla G_{B_{t}[j]}(x_{t - \tau_{t}^{s}}^{s}) )- \nabla G(x_{t - \tau_{t}^{s}}^{s})\|^2 
\nonumber\\
&&+ 2 \mathbb{E} \|\nabla G(x_{t - \tau_{t}^{s}}^{s})\|^2
\nonumber\\
&\overset{(\ref{young's})}{\leq}& \frac{2}{b} \sum_{j=1}^{b} \mathbb{E} \|\nabla G_{B_t[j]}(\widetilde{x}^s) - \nabla G_{B_{t}[j]}(x_{t - \tau_{t}^{s}}^{s}) - \nabla G(\widetilde{x}^s) + \nabla G(x_{t - \tau_{t}^{s}}^{s})\|^2 
+ 2B_{G}^2
\nonumber\\
&\overset{(\ref{expectation})}{=}& \frac{2}{b} \sum_{j=1}^{b} \mathbb{E} (\|\nabla G_{B_t[j]}(\widetilde{x}^s) - \nabla G_{B_{t}[j]}(x_{t - \tau_{t}^{s}}^{s})\|^2 - \|\nabla G(\widetilde{x}^s) - \nabla G(x_{t - \tau_{t}^{s}}^{s}))\|^2) 
\nonumber\\
&&+ 2B_{G}^2
\nonumber\\
&\leq& \frac{2}{b} \sum_{j=1}^{b} \mathbb{E} \|\nabla G_{B_t[j]}(\widetilde{x}^s) - \nabla G_{B_{t}[j]}(x_{t - \tau_{t}^{s}}^{s})\|^2 + 2B_{G}^2
\nonumber\\
&\overset{(\ref{young's})}{\leq}& \frac{2}{b} \sum_{j=1}^{b} \mathbb{E} 2 (\|\nabla G_{B_t[j]}(\widetilde{x}^s)\|^2 + \|\nabla G_{B_{t}[j]}(x_{t - \tau_{t}^{s}}^{s})\|^2) + 2B_{G}^2
\nonumber\\
&\leq& 8 B_{G}^2 + 2B_{G}^2
\nonumber\\
&=& 10 B_{G}^2,
\end{eqnarray}
where the last inequality comes from Assumption \ref{assumption:bounded gradient}.\\

\noindent \textbf{Proof for Lemma~\ref{lemma:4}.}

\noindent First, we bound $\mathbb{E} \|\nabla \widehat{f}(x_{t - \tau_{t}^{s}}^{s}) - \nabla \widehat{f}(x_{t}^{s})\|^2$ by
\begin{eqnarray} \label{shared_c_v_1}
&&\mathbb{E}\|\nabla \widehat{f}(x_{t - \tau_{t}^{s}}^{s}) - \nabla \widehat{f}(x_{t}^{s})\|^2 
\nonumber\\
&=& \mathbb{E}\|(\nabla \widehat{G}_{t - \tau_{t}^{s}}^{s})^{T} \nabla F_{i_t}(\widehat{G}_{t - \tau_{t}^{s}}^{s}) - (\nabla \widehat{G}_{t}^{s})^{T} \nabla F_{i_t}(\widehat{G}_{t}^{s}) \|^2
\nonumber\\
&=& \mathbb{E}\|(\nabla \widehat{G}_{t - \tau_{t}^{s}}^{s})^{T} \nabla F_{i_t}(\widehat{G}_{t - \tau_{t}^{s}}^{s}) - (\nabla \widehat{G}_{t - \tau_{t}^{s}}^{s})^{T} \nabla F_{i_t}(\widehat{G}_{t}^{s}) 
\nonumber\\
& &+ (\nabla \widehat{G}_{t - \tau_{t}^{s}}^{s})^{T} \nabla F_{i_t}(\widehat{G}_{t}^{s}) - (\nabla \widehat{G}_{t}^{s})^{T} \nabla F_{i_t}(\widehat{G}_{t}^{s}) \|^2
\nonumber\\
&\overset{(\ref{young's})(\ref{cauthy})}{\leq}& 2 \mathbb{E}\|\nabla \widehat{G}_{t - \tau_{t}^{s}}^{s}\|^2 \mathbb{E} \|\nabla F_{i_t}(\widehat{G}_{t - \tau_{t}^{s}}^{s}) - \nabla F_{i_t}(\widehat{G}_{t}^{s})\|^2
\nonumber\\
& & + 2\mathbb{E} \|\nabla  \widehat{G}_{t - \tau_{t}^{s}}^{s} - \widehat{G}_{t}^{s}\|^2 \mathbb{E} \|\nabla F_{i_t}(\widehat{G}_{t}^{s})\|^2.
\nonumber\\
&\leq& 20B_{G}^2 L_{F}^2 \mathbb{E} \|\widehat{G}_{t - \tau_{t}^{s}}^{s} - \widehat{G}_{t}^{s}\|^2 + 2B_{F}^2 \mathbb{E} \|\nabla \widehat{G}_{t - \tau_{t}^{s}}^{s} - \nabla \widehat{G}_{t}^{s}\|^2
\nonumber\\
&\overset{(\ref{young's})}{\leq}& 20 B_{G}^2 L_{F}^2 \frac{1}{b} \sum_{j=1}^{b} \mathbb{E} \|G_{B_{t}[j]}(x_{t}^{s})- G_{B_{t}[j]}^{s}(x_{t - \tau_{t}^{s}}^{s})\|^2
\nonumber\\
& & + \frac{2 B_{F}^2}{b} \sum_{j=1}^{b} \mathbb{E} \|\nabla G_{B_{t}[j]}(x_{t}^{s})- \nabla G_{B_{t}[j]}(x_{t - \tau_{t}^{s}}^{s})\|^2
\nonumber\\
&\leq& [20B_{G}^4 L_{F}^2 + 2B_{F}^2 L_{G}^2] \mathbb{E} \|x_{t}^{s} - x_{t - \tau_{t}^{s}}^{s}\|^2
\nonumber\\
&\overset{(\ref{young's})}{\leq}& [20B_{G}^4 L_{F}^2 + 2B_{F}^2 L_{G}^2] \cdot \frac{T \eta^2}{d_1} \sum_{j \in J(t)} \mathbb{E} \|\nabla \widehat{f}(x_{t-\tau_{t}^{s}}^s)\|^2,
\end{eqnarray}
where the second inequality follows Lemma~\ref{lemma:3} and the forth inequality comes from Assumption \ref{assumption:bounded gradient} and Assumption \ref{assumption:smoothness}. It follows that
\begin{eqnarray}
&&\mathbb{E} \|\nabla \widehat{f}(x_{t - \tau_{t}^{s}}^{s})\|^2 
\nonumber\\
&\leq& 2 \mathbb{E} \|\nabla \widehat{f}(x_{t - \tau_{t}^{s}}^{s}) - \nabla \widehat{f}(x_{t}^{s})\|^2 + 2 \mathbb{E} \|\nabla \widehat{f}(x_{t}^{s})\|^2
\nonumber\\
&\overset{(\ref{shared_c_v_1})}{\leq}& [40B_{G}^4 L_{F}^2 + 4 B_{F}^2 L_{G}^2] \frac{T \eta^2}{d_1} \sum_{j \in J(t)} \mathbb{E} \|\nabla \widehat{f}(x_{t- \tau_{t}^{s}}^{s})\|^2 + 2 \mathbb{E} \|\nabla \widehat{f}(x_{t}^{s})\|^2.
\end{eqnarray}
Summing up this inequality from $t=0$ to $t=K-1$ yields
\begin{eqnarray}
&&\sum_{t=0}^{K-1} \mathbb{E} \|\nabla \widehat{f}(x_{t - \tau_{t}^{s}}^{s})\|^2 
\nonumber\\
&\leq& [40B_{G}^4 L_{F}^2 + 4 B_{F}^2 L_{G}^2] \frac{T \eta^2}{d_1} \sum_{t=0}^{K-1} \sum_{j \in J(t)} \mathbb{E} \|\nabla \widehat{f}(x_{t - \tau_{t}^{s}}^{s})\|^2 + 2 \sum_{t=0}^{K-1} \mathbb{E} \|\nabla \widehat{f}(x_{t}^{s})\|^2
\nonumber\\
&\leq& [40B_{G}^4 L_{F}^2 + 4 B_{F}^2 L_{G}^2] \frac{T^2 \eta^2}{d_1} \sum_{t=0}^{K-1} \mathbb{E} \|\nabla \widehat{f}(x_{t - \tau_{t}^{s}}^{s})\|^2 + 2 \sum_{t=0}^{K-1} \mathbb{E} \|\nabla \widehat{f}(x_{t}^{s})\|^2,
\end{eqnarray}
where the last inequality can be obtained by using a simple augment and the time delays are at most $T$. Then, $\sum_{t=0}^{K-1} \mathbb{E} \|\nabla \widehat{f}(x_{t - \tau_{t}^{s}}^{s})\|^2$ can be bounded
\begin{eqnarray}
\sum_{t=0}^{K-1} \mathbb{E} \|\nabla \widehat{f}(x_{t - \tau_{t}^{s}}^{s})\|^2 \leq \frac{2}{1- [40B_{G}^4 L_{F}^2 + 4 B_{F}^2 L_{G}^2] \frac{T^2 \eta^2}{d_1} } \sum_{t=0}^{K-1} \mathbb{E} \|\nabla \widehat{f}(x_{t}^{s})\|^2.
\end{eqnarray}\\

\noindent \textbf{Proof of Theorem~\ref{theorem_shared_1:1}.}

\noindent For AsyVRSC-Shared, the iteration at time $t$ is $(x_{t+1}^{s})_{k_{t}} = (x_{t}^{s})_{k_{t}} - \eta (\nabla \widehat{f}(x_{t - \tau_{t}^{s}}^{s}))_{k_{t}}$. Since we assume the time delays have an upper bound $T$, $x_{t-\tau_{t}^{s}}^{s}$ can be expressed as  $x_{t-\tau_{t}^{s}}^{s} = x_{t}^{s} - \sum_{j \in J(t)}( x_{j+1}^{s} - x_{j}^{s})$, where $J(t) \in \{ t-1, t-2, ..., t-T \}$.
We start by decomposing the expectation of $\| x_{t+1}^{s} - x^{*} \|^2$ as
\begin{eqnarray} \label{shared_first}
&&\mathbb{E} \|x_{t+1}^{s} - x^{*}\|^{2} 
\nonumber\\
&=& \mathbb{E} \|x_{t+1}^{s} - x_{t}^{s} + x_{t}^{s} - x^{*}\|^2
\nonumber\\
&=& \mathbb{E} \|x_{t+1}^{s} - x_{t}^{s}\|^2 + \mathbb{E} \|x_{t}^{s} - x^{*}\|^2 + 2\mathbb{E} \langle x_{t+1}^{s} - x_{t}^{s},  x_{t}^{s}- x^{*} \rangle
\nonumber\\
&=& \frac{\eta^2}{d_1} \mathbb{E} \|\nabla \widehat{f}(x_{t - \tau_{t}^{s}}^{s})\|^2 + \mathbb{E} \|x_{t}^{s} - x^{*}\|^2 + \frac{2 \eta}{d_1} \underbrace{\mathbb{E} \langle \nabla \widehat{f}(x_{t - \tau_{t}^{s}}^{s}), x^{*} - x_{t}^{s} \rangle }_{T_{1}}.
\end{eqnarray}
We then bound $T_{1}$ by
\begin{eqnarray} \label{shared_T1_1}
T_{1} &=& \mathbb{E} \langle x^{*} - x_{t}^{s}, (\nabla \widehat{G}_{t - \tau_{t}^{s}}^{s})^{T} \nabla F_{i_{t}}(\widehat{G}_{t - \tau_{t}^{s}}^{s}) \rangle
\nonumber\\
&=& \underbrace{\mathbb{E} \langle x^{*} - x_{t}^{s}, (\nabla \widehat{G}_{t - \tau_{t}^{s}}^{s})^{T} \nabla F_{i_{t}}(\widehat{G}_{t - \tau_{t}^{s}}^{s}) - \nabla f(x_{t - \tau_{t}^{s}}^{s}) \rangle }_{T_{2}} 
\nonumber\\
&&+ \underbrace{\mathbb{E} \langle x^{*} - x_{t}^{s},  \nabla f(x_{t - \tau_{t}^{s}}^{s}) \rangle }_{T_{3}}.
\end{eqnarray}
We proceed to bound $T_2$ by
\begin{eqnarray}\label{shared_T2_1}
T_{2} 
&=& \mathbb{E} \langle x^{*} - x_{t}^{s}, (\nabla G(x_{t - \tau_{t}^{s}}^{s}))^{T} \nabla F_{i_{t}}(\widehat{G}_{t - \tau_{t}^{s}}^{s}) - (\nabla G(x_{t - \tau_{t}^{s}}^{s}))^{T} \nabla F_{i_{t}}(G_{t - \tau_{t}^{s}}^{s}) \rangle 
\nonumber\\
&\overset{(\ref{young's2})}{\leq}& \frac{\alpha}{2} \mathbb{E} \|x_{t}^{s} - x^{*}\|^2 + \frac{1}{2\alpha} \mathbb{E}\|(\nabla G(x_{t - \tau_{t}^{s}}^{s}))^{T} \nabla F_{i_{t}}(\widehat{G}_{t - \tau_{t}^{s}}^{s}) 
\nonumber\\
&&- (\nabla G(x_{t - \tau_{t}^{s}}^{s}))^{T} \nabla F_{i_{t}}(G_{t - \tau_{t}^{s}}^{s})\|^2
\nonumber\\
&\overset{(\ref{cauthy})}{\leq}&  \frac{\alpha}{2} \mathbb{E} \|x_{t}^{s} - x^{*}\|^2 + \frac{1}{2\alpha} \mathbb{E} \|\nabla G(x_{t - \tau_{t}^{s}}^{s})\|^2 \|\nabla F_{i_{t}}(\widehat{G}_{t - \tau_{t}^{s}}^{s}) - \nabla F_{i_{t}}(G_{t - \tau_{t}^{s}}^{s})\|^2
\nonumber\\
&\leq& \frac{\alpha}{2} \mathbb{E} \|x_{t}^{s} - x^{*}\|^2 + \frac{B_{G}^2 L_{F}^2}{2\alpha} \mathbb{E} \|\widehat{G}_{t - \tau_{t}^{s}}^{s} - G_{t - \tau_{t}^{s}}^{s}\|^2,
\nonumber\\
&\overset{(\ref{lemma:bound_G})}{\leq}& ( \frac{\alpha}{2} + \frac{3 B_{G}^4 L_{F}^2}{2a\alpha} ) \mathbb{E} \|x_{t}^{s} - x^{*}\|^2 + \frac{3 B_{G}^4 L_{F}^2}{2a\alpha} \mathbb{E} (\|\widetilde{x}^{s} - x^{*}\|^2 + \|x_{t}^{s} - x_{t - \tau_{t}^{s}}^{s}\|^2). \hspace{15pt}
\end{eqnarray}
where the third inequality comes from Assumption \ref{assumption:bounded gradient} and Assumption \ref{assumption:smoothness} and the last inequality follows Lemma~\ref{lemma:2}. Also, we bound $T_{3}$ by
\begin{eqnarray} \label{shared_T3_1}
- T_{3} &=& \mathbb{E} \langle x_{t} - x^{*}, \nabla f(x_{t - \tau_{t}^{s}}^{s}) \rangle
\nonumber\\
&=& \mathbb{E} \langle x_{t - \tau_{t}^{s}}^{s} - x^{*}, \nabla f(x_{t - \tau_{t}^{s}}^{s}) \rangle + \mathbb{E} \langle x_{t}^{s} - x_{t - \tau_{t}^{s}}^{s}, \nabla f(x_{t - \tau_{t}^{s}}^{s}) \rangle
\nonumber\\
&\geq& \mathbb{E}(f(x_{t - \tau_{t}^{s}}^{s}) - f(x^{*})) + \mathbb{E} \langle x_{t}^{s} - x_{t - \tau_{t}^{s}}^{s}, \nabla f(x_{t - \tau_{t}^{s}}^{s}) \rangle
\nonumber\\
&\geq& \mathbb{E}(f(x_{t - \tau_{t}^{s}}^{s}) - f(x^{*})) + \mathbb{E} (f(x_{t}^{s}) - f(x_{t - \tau_{t}^{s}}^{s})) - \frac{L_{f}}{2} \mathbb{E} \|x_{t}^{s} - x_{t - \tau_{t}^{s}}^{s}\|^2
\nonumber\\
&=& \mathbb{E}(f(x_{t}^{s}) - f(x^{*})) - \frac{L_{f}}{2} \mathbb{E} \|x_{t}^{s} - x_{t - \tau_{t}^{s}}^{s}\|^2
\nonumber\\
&=& \mathbb{E}(f(x_{t}^{s}) - f(x^{*})) - \frac{L_{f}}{2} \mathbb{E} \|\sum_{j \in J(t)} ( x_{j+1}^{s} - x_{j}^{s})\|^2
\nonumber\\
&\overset{(\ref{young's})}{\geq}& \mathbb{E}(f(x_{t}^{s}) - f(x^{*})) - \frac{L_{f} T}{2} \sum_{j \in J(t)} \mathbb{E} \|x_{j+1}^{s} - x_{j}^{s}\|^2
\nonumber\\
&=& \mathbb{E}( f(x_{t}^{s}) - f(x^{*})) - \frac{L_{f} \eta^2 T}{2 d_1} \sum_{j \in J(t)} \mathbb{E} \|\nabla \widehat{f}(x_{j - \tau_{j}^{s}}^{s})\|^2,
\end{eqnarray}
where the first and the second inequalities come from the smoothness and convexity of $f(x)$, respectively. Putting (\ref{shared_T2_1}) and (\ref{shared_T3_1}) back to (\ref{shared_T1_1}), we have
\begin{eqnarray} \label{shared_T1_2}
T_{1} &\leq& (\frac{\alpha}{2} + \frac{3 B_{G}^4 L_{F}^2}{2a \alpha}) \mathbb{E} \|x_{t}^{s} - x^{*}\|^2 + \frac{3 B_{G}^4 L_{F}^2}{2a\alpha} \mathbb{E} (\|\widetilde{x}_{t}^{s} - x^{*}\|^2 + \|x_{t}^{s} - x_{t - \tau_{t}^{s}}^{s}\|^2 )
\nonumber\\
& & +\frac{L_{f} \eta^2 T}{2d_1} \sum_{j \in J(t)} \mathbb{E} \|\nabla \widehat{f}(x_{j - \tau_{j}^{s}}^{s})\|^2 - \mathbb{E} (f(x_{t}^{s}) - f(x^{*})).
\end{eqnarray}

Applying the upper bound of $T_1$ in (\ref{shared_T1_2}) to (\ref{shared_first}) yields
\begin{eqnarray}
&&\mathbb{E} \|x_{t+1}^{s} - x^*\|^2 
\nonumber\\
&\leq& \mathbb{E} \|x_{t}^{s} - x^{*}\|^2 + \frac{\eta^2}{d_1} \mathbb{E} \|\nabla \widehat{f}(x_{t - \tau_{t}^{s}}^{s})\|^2 + \frac{1}{d_1} (\alpha \eta + \frac{3 \eta B_{G}^4 L_{F}^2}{a \alpha}) \mathbb{E} \|x_{t}^{s} - x^{*}\|^2
\nonumber\\
&& +\frac{3 \eta B_{G}^4 L_{F}^2} {d_1 a\alpha} \mathbb{E} (\|\widetilde{x}_{t}^{s} - x^{*}\|^2 + \|x_{t}^s - x_{t - \tau_{t}^{s}}^{s}\|^2) + \frac{L_{f} \eta^3 T}{d_1^2} \sum_{j \in J(t)} \mathbb{E} \|\nabla \widehat{f}(x_{j - \tau_{j}^{s}}^{s})\|^2
\nonumber\\
&&- \frac{2 \eta}{d_1} \mathbb{E} (f(x_{t}^{s}) - f(x^{*})),
\nonumber\\
&\overset{(\ref{young's})}{\leq}& \mathbb{E} \|x_{t}^{s} - x^{*}\|^2 + \frac{\eta^2}{d_1} \mathbb{E} \|\nabla \widehat{f}(x_{t - \tau_{t}^{s}}^{s})\|^2 + \frac{1}{d_1} (\alpha \eta + \frac{3 \eta B_{G}^4 L_{F}^2}{a \alpha}) \mathbb{E} \|x_{t}^{s} - x^{*}\|^2
\nonumber\\
&& +\frac{3 \eta B_{G}^4 L_{F}^2} {d_1 a\alpha} \mathbb{E} \|\widetilde{x}_{t}^{s} - x^{*}\|^2  + (\frac{L_{f} \eta^3 T}{d_1^2} + \frac{3 \eta^3 B_{G}^4 L_{F}^2 T} {d_1^2 a\alpha} )\sum_{j \in J(t)} \mathbb{E} \|\nabla \widehat{f}(x_{j - \tau_{j}^{s}}^{s})\|^2
\nonumber\\
&&- \frac{2 \eta}{d_1} \mathbb{E} (f(x_{t}^{s}) - f(x^{*})).
\end{eqnarray}
Summing up this inequality from $t=0$ to $t=K-1$, we obtain
\begin{eqnarray} \label{shared_first_2}
&&\mathbb{E} \|x_{K}^{s} - x^{*}\|^2 
\nonumber\\
&\leq& \mathbb{E} \|\widetilde{x}^{s} - x^{*}\|^2 + \frac{1}{d_1}(\alpha \eta + \frac{3 \eta B_{G}^4 L_{F}^2}{a \alpha}) \sum_{t=0}^{K-1}\mathbb{E} \|x_{t}^{s} - x^{*}\|^2 
\nonumber\\
&& + \frac{\eta^2}{d_1} \sum_{t=0}^{K-1} \mathbb{E} \| \nabla \widehat{f}(x_{t - \tau_{t}^{s}}^{s})\|^2 + \frac{3 \eta B_{G}^4 L_{F}^2}{d_1 a \alpha} K \mathbb{E} \|\widetilde{x}^s - x^{*}\|^2 
\nonumber\\
&&+ (\frac{L_{f} \eta^3 T}{d_1^2} + \frac{3 \eta^3 B_{G}^4 L_{F}^2 T}{d_1^2 a \alpha}) \sum_{t=0}^{K-1} \sum_{j \in J(t)} \mathbb{E} \|\nabla \widehat{f}(x_{t - \tau_{j}^{s}}^{s})\|^2
\nonumber\\
&&- \frac{2 \eta }{d_1} \sum_{t=0}^{K-1} \mathbb{E} (f(x_{t}^{s}) - f(x^{*}))
\nonumber\\
&\leq& (1 + \frac{3 \eta B_{G}^4 L_{F}^2 K}{d_1 a \alpha}) \mathbb{E} \|\widetilde{x}^{s} - x^{*}\|^2  + \frac{1}{d_1}(\alpha T + \frac{3 \eta B_{G}^4 L_{F}^2}{a \alpha}) \sum_{t=0}^{K-1}\mathbb{E} \|x_{t}^{s} - x^{*}\|^2
\nonumber\\
&& +\frac{\eta^2}{d_1} \sum_{t=0}^{K-1} \mathbb{E} \| \nabla \widehat{f}(x_{t - \tau_{t}^{s}}^{s})\|^2 + (\frac{L_{f} \eta^3 \tau^2}{d_1^2} + \frac{3 \eta^3 B_{G}^4 L_{F}^2 \tau^2}{d_1^2 a \alpha}) \sum_{t=0}^{K-1} \mathbb{E} \|\nabla \widehat{f}(x_{t - \tau_{t}^{s}}^{s})\|^2
\nonumber\\
&&- \frac{2 \eta }{d_1} \sum_{t=0}^{K-1} \mathbb{E} (f(x_{t}^{s}) - f(x^{*}))
\nonumber\\
&\leq& (1 + \frac{3 \eta B_{G}^4 L_{F}^2 K}{d_1 a \alpha}) \frac{\mu_{f}}{2} \mathbb{E} (f(\widetilde{x}^s) - f(x^{*}) ) 
\nonumber\\
&& - (\frac{2 \eta }{d_1} - \frac{1}{d_1}(\alpha \eta + \frac{3 \eta B_{G}^4 L_{F}^2}{a \alpha})\frac{\mu_{f}}{2} ) \mathbb{E} \sum_{t=0}^{K-1} ( f(x_{t}^{s}) - f(x^{*}) )
\nonumber\\
&& + (\frac{\eta^2}{d_1} + \frac{L_{f} \eta^3 T^2}{d_1^2} + \frac{3 \eta^3 B_{G}^4 L_{F}^2 T^2}{d_1^2 a \alpha}) \sum_{t=0}^{K-1} \mathbb{E} \|\nabla \widehat{f}(x_{t - \tau_{t}^{s}}^{s})\|^2. 
\end{eqnarray}
where the second inequality can be obtained by using a simple counting augment and the fact that the time delays are at most $T$, and the third inequality follows from that $f(x)$ is $\mu_{f}$-strongly convex. 
Combining with Lemma~\ref{lemma:1} and Lemma~\ref{lemma:4} and setting $\alpha = \frac{\mu_{f}}{8}$, we obtain
\begin{eqnarray}
&&\mathbb{E} \|x_{K}^{s} - x^{*}\|^2
\nonumber\\
&\leq& \frac{2}{\mu_{f}} (1 + \frac{24 \eta B_{G}^4 L_{F}^2 K}{d_1 a \mu_f}) \mathbb{E} (f(\widetilde{x}^s) - f(x^{*})) 
\nonumber\\
&&- (\frac{2 \eta }{d_1} - \frac{1}{d_1}( \frac{\eta \mu_f}{8} + \frac{24 \eta B_{G}^4 L_{F}^2}{a \mu_f})\frac{\mu_{f}}{2} ) \mathbb{E} \sum_{t=0}^{K-1} ( f(x_{t}^{s}) - f(x^{*}) )
\nonumber\\
&& + (\frac{\eta^2}{d_1} + \frac{L_{f} \eta^3 T^2}{d_1^2} + \frac{24 \eta^3 B_{G}^4 L_{F}^2 T^{2}}{d_1^2 a \mu_f}) \frac{2}{1- [40B_{G}^4 L_{F}^2 + 4 B_{F}^2 L_{G}^2] \frac{T^2 \eta^2}{d_1} }
\nonumber\\
&& (\frac{64}{\mu_{f}}(\frac{B_{G}^4 L_{F}^2}{a} + \frac{B_{F}^{2} L_{G}^{2}}{b} ) + 8 L_{f}) \sum_{t=0}^{K-1} \mathbb{E}(f(x_{t}^{s}) - f(x^{*}) + f(\widetilde{x}^{s}) - f(x^{*})).
\end{eqnarray}
Discarding the left hand side and setting $U, P, Q, R$ as (\ref{theorem_shared_1:parameter}),
we have
\begin{equation} \label{shared_last_1}
\frac{1}{K} \sum_{t=0}^{K-1} \mathbb{E} ( f(x_{t}^{s}) - f(x^{*}) ) \leq \frac{\frac{\mu_{f}}{2}  + PQRK + U}{\frac{7 \eta K}{4d_1} - PQRK - U} \mathbb{E}(f(\widetilde{x}^s) - f(x^{*})).
\end{equation}
Since $f(x)$ is convex and $\widetilde{x}^{s+1} = \mathbb{E}_{t \in \{ 0,\cdots,K-1\}} ~x_{t}^{s}$, we can bound the left hand side of (\ref{shared_last_1}) by using Jensen's inequality
\begin{equation} \label{jeason}
\frac{1}{K} \sum_{t=0}^{K-1} ( f(x_{t}^{s}) - f(x^{*}) ) \geq f(\widetilde{x}^{s+1}) - f(x^{*}).
\end{equation}
Substituting (\ref{jeason}) to (\ref{shared_last_1}), we complete the proof and obtain
\begin{equation} \label{shared_last_2}
\mathbb{E} ( f(\widetilde{x}^{s+1}) - f(x^{*}) ) \leq \frac{\frac{\mu_{f}}{2}  + PQRK + U}{\frac{7 \eta K}{4d_1} - PQRK - U} \mathbb{E}(f(\widetilde{x}^s) - f(x^{*})).
\end{equation}

\noindent \textbf{Proof of Corollary~\ref{corollary_shared_1:1}.}

\noindent The main idea to choose the parameters in Theorem \ref{theorem_shared_1:1} is to ensure the geometric convergence parameter $\frac{\frac{\mu_{f}}{2}  + PQRK + U}{\frac{7 \eta K}{4d_1} - PQRK - U} < 1$. If $T$ can be bounded by $\sqrt{d_1}$ and by choosing 
$a = \max\{ \frac{1024 B_{G}^4 L_{F}^2}{\mu_{f}^2}, \frac{32B_{G}^4 L_F^2}{5\mu_f L_f}\}, b = \frac{32 B_{F}^2 L_{G}^2}{\mu_{f} L_{f}} $,
we can bound $U$, $P$, $Q$ and $R$ as
\begin{equation}
U = \frac{48 \eta B_{G}^{4} L_{F}^{2} K}{d a \mu_{f}^{2}} \leq \frac{3\eta K}{64 d_1},
\end{equation}
\begin{equation}
P = \frac{\eta^2}{d_1} + \frac{L_{f} \eta^3 T^2}{d_1^2} + \frac{24 \eta^3 B_{G}^4 L_{F}^2 T^{2}}{d_1^2 a \mu_{f}}
\leq \frac{\eta^2}{d_1} + \frac{L_f \eta^3}{d_1} + \frac{15 L_f \eta^3}{4d_1}
= \frac{\eta^2}{d_1} + \frac{19 L_f \eta^3}{4d_1},
\end{equation}
\begin{equation}
Q = \frac{2}{1- [40B_{G}^4 L_{F}^2 + 4 B_{F}^2 L_{G}^2] \frac{T^2 \eta^2}{d_1} } \leq \frac{2}{1- [40B_{G}^4 L_{F}^2 + 4 B_{F}^2 L_{G}^2] \eta^2 },
\end{equation}
\begin{equation}
R = \frac{64}{\mu_{f}}(\frac{B_{G}^4 L_{F}^2}{a} + \frac{B_{F}^{2} L_{G}^{2}}{b} ) + 8 L_{f} \leq \frac{64}{\mu_f} (\frac{5 \mu_f L_f}{32} + \frac{\mu_f L_f}{ 32}) + 8L_f = 20 L_f.
\end{equation}
We choose $\eta$ as
\begin{equation}
\eta = \min \{ \frac{1}{9 B_{G}^2 L_{F}}, \frac{1}{9B_{F} L_{G}},\frac{1}{320 L_{f}}\}.
\end{equation}
Then $U$, $P$, $Q$ and $R$ can be further bounded
\begin{equation}
U \leq \frac{3K}{20480 L_f d_1}, P \leq \frac{1299}{131072000  L_f^2 d_1}, Q \leq \frac{162}{37}, R \leq 20 L_f,
\end{equation}
\begin{equation}
PQR \leq \frac{105219}{121241600 L_f d_1}.
\end{equation}
At last, choosing $K = \frac{1024 L_{f} d_1}{\mu_{f}}$, we have
\begin{equation}
\frac{\frac{\mu_{f}}{2}  + PQRK + U}{\frac{7 \eta K}{4d_1} - PQRK - U} \leq \frac{\mu_f (\frac{1}{2} + \frac{105219}{118400} + \frac{3}{20})}{ \mu_{f} (\frac{28}{5} - \frac{105219}{118400} + \frac{3}{20})} \approx 0.666182 \leq \frac{2}{3}.
\end{equation}
We can obtain a linear convergence rate $\frac{2}{3}$. This completes the proof.

\section*{Appendix C: Convergence Analysis for Section 4.2}
\textbf{Proof of Lemma~\ref{lemma:5}.}

\noindent First, we bound $\mathbb{E}\|\nabla \widehat{f}(x_{t - \tau_{t}^{s}}^{s}) - \nabla \widehat{f}(x_{t}^{s})\|^2$ as
\begin{eqnarray} \label{distribute_c_v_1}
&&\mathbb{E}\|\nabla \widehat{f}(x_{t - \tau_{t}^{s}}^{s}) - \nabla \widehat{f}(x_{t}^{s})\|^2 
\nonumber\\
&=& \mathbb{E}\|(\nabla \widehat{G}_{t - \tau_{t}^{s}}^{s})^{T} \nabla F_{i_t}(\widehat{G}_{t - \tau_{t}^{s}}^{s}) - (\nabla \widehat{G}_{t}^{s})^{T} \nabla F_{i_t}(\widehat{G}_{t}^{s}) \|^2
\nonumber\\
&=& \mathbb{E}\|(\nabla \widehat{G}_{t - \tau_{t}^{s}}^{s})^{T} \nabla F_{i_t}(\widehat{G}_{t - \tau_{t}^{s}}^{s}) - (\nabla \widehat{G}_{t - \tau_{t}^{s}}^{s})^{T} \nabla F_{i_t}(\widehat{G}_{t}^{s}) 
\nonumber\\
& &+ (\nabla \widehat{G}_{t - \tau_{t}^{s}}^{s})^{T} \nabla F_{i_t}(\widehat{G}_{t}^{s}) - (\nabla \widehat{G}_{t}^{s})^{T} \nabla F_{i_t}(\widehat{G}_{t}^{s}) \|^2
\nonumber\\
&\overset{(\ref{young's})(\ref{cauthy})}{\leq}& 2 \mathbb{E}\|\nabla \widehat{G}_{t - \tau_{t}^{s}}^{s}\|_{i_t}^2 \mathbb{E} \|\nabla F_{i_t}(\widehat{G}_{t - \tau_{t}^{s}}^{s}) - \nabla F_{i_t}(\widehat{G}_{t}^{s})\|_{i_t}^2
\nonumber\\
& & + 2\mathbb{E} \|\nabla  \widehat{G}_{t - \tau_{t}^{s}}^{s} - \nabla \widehat{G}_{t}^{s}\|_{i_t}^2 \mathbb{E} \|\nabla F_{i_t}(\widehat{G}_{t}^{s})\|_{i_t}^2
\nonumber\\
&\leq& 20B_{G}^{2} \Delta \mathbb{E} \|\nabla F_{i_t}(\widehat{G}_{t - \tau_{t}^{s}}^{s}) - \nabla F_{i_t}(\widehat{G}_{t}^{s})\|^2 + 2 \Delta B_{F}^{2} \mathbb{E} \|\nabla  \widehat{G}_{t - \tau_{t}^{s}}^{s} - \nabla \widehat{G}_{t}^{s}\|^2
\nonumber\\
&\leq& 20B_{G}^{2} \Delta L_{F}^{2} \mathbb{E} \|\widehat{G}_{t - \tau_{t}^{s}}^{s} - \widehat{G}_{t}^{s}\|^2 + 2 \Delta B_{F}^{2} \mathbb{E} \|\nabla  \widehat{G}_{t - \tau_{t}^{s}}^{s} - \nabla \widehat{G}_{t}^{s}\|^2
\nonumber\\
&\leq& (20 B_{G}^{4} L_{F}^{2} \Delta + 2 B_{F}^{2} L_{G}^{2} \Delta) T \eta^{2} \sum_{l=t-\tau_{t}^{s}}^{t-1} \mathbb{E} \| \nabla \widehat{f}(x_{l}^{s}) \|^2,
\end{eqnarray}
where $\| \cdot\|_{i_{t}}$ denotes the support of $\nabla F_{i_{t}}(x)$, and the second inequality comes from Lemma~\ref{lemma:3} and Assumption~\ref{assumption:sparsity}. The last inequality is similar to (\ref{shared_c_v_1}).
Following (\ref{distribute_c_v_1}), we can bound $\mathbb{E} \|\nabla \widehat{f}(x_{t - \tau_{t}^{s}}^{s})\|^2$ by
\begin{eqnarray} \label{bound_d}
&&\mathbb{E} \|\nabla \widehat{f}(x_{t - \tau_{t}^{s}}^{s})\|^2 
\nonumber\\
& \overset{(\ref{young's})}{\leq}& 2 \mathbb{E} \|\nabla \widehat{f}(x_{t - \tau_{t}^{s}}^{s}) - \nabla \widehat{f}(x_{t}^{s})\|^2 + 2 \mathbb{E} \|\nabla \widehat{f}(x_{t}^{s})\|^2
\nonumber\\
&\leq& [40B_{G}^4 L_{F}^2 + 4B_{F}^2 L_{G}^2]\Delta \eta^2 T \sum_{l=t - \tau_{t}^{s}}^{t-1} \mathbb{E} \| \nabla \widehat{f}(x_{l}^{s}) \|^2 + 2 \mathbb{E} \|\nabla \widehat{f}(x_{t}^{s})\|^2.
\end{eqnarray}
Summing up this inequality from $t=0$ to $t=K-1$, we get
\begin{eqnarray}\label{summary_bound_d}
&&\sum_{t=0}^{K-1} \mathbb{E} \|\nabla \widehat{f}(x_{t - \tau_{t}^{s}}^{s})\|^2  
\nonumber\\
&\leq& [40B_{G}^4 L_{F}^2 + 4B_{F}^2 L_{G}^2]\Delta \eta^2 T \sum_{t=0}^{K-1} \sum_{l=t-\tau_{t}^{s}}^{t-1} \mathbb{E} \|\nabla \widehat{f}(x_{l}^{s})\|^2 +2 \sum_{t=0}^{K-1} \mathbb{E} \|\nabla \widehat{f}(x_{t}^{s})\|^2
\nonumber\\
&\leq& [40B_{G}^4 L_{F}^2 + 4B_{F}^2 L_{G}^2]\Delta \eta^2 T^2 \sum_{t=0}^{K-1} \mathbb{E} \|\nabla \widehat{f}(x_{t-\tau_{t}^{s}}^{s})\|^2 +2 \sum_{t=0}^{K-1} \mathbb{E} \|\nabla \widehat{f}(x_{t}^{s})\|^2.
\end{eqnarray}
Then, $\sum_{t=0}^{K-1} \mathbb{E} \|\nabla \widehat{f}(x_{t - \tau_{t}^{s}}^{s})\|^2$ can be bounded by
\begin{eqnarray} \label{shared_T5_1}
\sum_{t=0}^{K-1} \mathbb{E} \|\nabla \widehat{f}(x_{t - \tau_{t}^{s}}^{s})\|^2 \leq \frac{2}{1- [40B_{G}^4 L_{F}^2 + 4 B_{F}^2 L_{G}^2] \Delta \eta^2 T^2} \sum_{t=0}^{K-1} \mathbb{E} \|\nabla \widehat{f}(x_{t}^{s})\|^2.
\end{eqnarray}\\\\

\noindent \textbf{Proof of Theorem~\ref{theorem_ditribute_1:1}.}

\noindent First, the iteration at time $t$ of epoch $s$  is
\begin{equation} \label{iter}
x_{t+1}^{s} = x_{t}^{s} - \eta \nabla \widehat{f}(x_{t- \tau_{t}^{s}}^{s}),
\end{equation}
and then we have
\begin{equation} \label{Theorem_2}
\mathbb{E}\|x_{t+1}^{s} - x^{*}\|^2 = \mathbb{E}\|x_{t}^{s} - x^{*}\|^{2} + \eta^2  \mathbb{E} \|\nabla \widehat{f}(x_{t - \tau_{t}^{s}}^{s})\|^2 + 2\eta \underbrace{ \mathbb{E} \langle x^{*} - x_{t}^{s},\nabla \widehat{f}(x_{t - \tau_{t}^{s}}^{s})\rangle}_{T_4}.
\end{equation}
We then bound $T_{4}$ by
\begin{eqnarray}\label{distribute_T7}
T_{4} & = & \mathbb{E} \langle x^{*} - x_{t}^{s}, (\nabla \widehat{G}_{t - \tau_{t}^{s}}^{s})^{T} \nabla F_{i_t}(\widehat{G}_{t - \tau_{t}^{s}}^{s}) \rangle
\nonumber\\
&=& \underbrace{\mathbb{E} \langle x^{*} - x_{t}^{s}, (\nabla \widehat{G}_{t - \tau_{t}^{s}}^{s})^{T} \nabla F_{i_t}(\widehat{G}_{t - \tau_{t}^{s}}^{s}) - \nabla f(x_{t - \tau_{t}^{s}}^{s}) \rangle }_{T_5} 
\nonumber\\
&&+ \underbrace{\mathbb{E} \langle x^{*} - x_{t}^{s}, \nabla f(x_{t - \tau_{t}^{s}}^{s}) \rangle }_{T_6}.
\end{eqnarray}
Following (\ref{young's2}) and (\ref{cauthy}), we can bound $T_5$ as
\begin{eqnarray} \label{distribute_T8}
&&-T_{5} 
\nonumber\\
&\overset{(\ref{young's2})}{\geq}& -\frac{\alpha}{2} \mathbb{E} \| x^{*} - x_{t}^{s} \|^2 - \frac{1}{2\alpha} \mathbb{E} \| (\nabla G_{j_{t}}(x_{t - \tau_{t}^{s}}^{s}))^{T} \nabla F_{i_{t}}(\widehat{G}_{t-\tau_{t}^{s}}^{s}) 
\nonumber\\
&&- (\nabla G_{j_{t}}(x_{t - \tau_{t}^{s}}^{s}))^{T} \nabla F_{i_{t}}(G(x_{t-\tau_{t}^{s}}^{s})) \|^{2}
\nonumber\\
&\overset{(\ref{cauthy})}{\geq}&  -\frac{\alpha}{2} \mathbb{E} \| x^{*} - x_{t}^{s} \|^2 - \frac{1}{2\alpha} \mathbb{E} \| \nabla G_{j_{t}}(x_{t - \tau_{t}^{s}}^{s}) \|_{i_{t}}^2 \| \nabla F_{i_{t}}(\widehat{G}_{t-\tau_{t}^{s}}^{s}) - \nabla F_{i_{t}}(G(x_{t-\tau_{t}^{s}}^{s})) \|_{i_{t}}^{2}
\nonumber\\
&\geq& -\frac{\alpha}{2} \mathbb{E} \| x^{*} - x_{t}^{s} \|^2 - \frac{\Delta B_{G}^{2}}{2\alpha} \mathbb{E} \| \nabla F_{i_{t}}(\widehat{G}_{t-\tau_{t}^{s}}^{s}) - \nabla F_{i_{t}}(G(x_{t-\tau_{t}^{s}}^{s})) \|^{2}
\nonumber\\
&\geq& -\frac{\alpha}{2} \mathbb{E} \| x^{*} - x_{t}^{s} \|^2 - \frac{\Delta B_{G}^{2} L_{F}^{2}}{2\alpha} \underbrace{\mathbb{E} \| \widehat{G}_{t-\tau_{t}^{s}}^{s} - G(x_{t-\tau_{t}^{s}}^{s}) \|^{2} }_{T_{7}},
\end{eqnarray}
where the third inequality comes from Assumption \ref{assumption:sparsity}. Using the bound given by (\ref{shared_T4_2}), we can bound $T_{7}$ with
\begin{eqnarray} \label{distribute_T10}
T_{7} &\leq& \frac{3 B_{G}^2}{a}\mathbb{E}(\|\widetilde{x}^s  - x^{*}\|^2 + \|x_{t}^{s} - x^{*}\|^2 + \| x_{t - \tau_{t}^{s}}^{s} - x_{t}^{s}\|^2).
\end{eqnarray}
Bring (\ref{distribute_T10}) into (\ref{distribute_T8}), we can get
\begin{eqnarray} \label{distribute_T8_2}
T_5 &\leq& \frac{\alpha}{2} \mathbb{E} \|x_{t}^{s} - x^{*}\|^2 + \frac{3B_{G}^4L_{F}^2 \Delta}{2a\alpha} \mathbb{E} (\|\widetilde{x}^s  - x^{*}\|^2 + \|x_{t}^{s} - x^{*}\|^2 + \| x_{t - \tau_{t}^{s}}^{s} - x_{t}^{s}\|^2)
\nonumber\\
&\leq& (\frac{\alpha}{2} + \frac{3B_{G}^4L_{F}^2 }{2a\alpha})\mathbb{E} \|x_{t}^{s} - x^{*}\|^2 + \frac{3B_{G}^4L_{F}^2 }{2a\alpha} \mathbb{E} \|\widetilde{x}^s  - x^{*}\|^2 
\nonumber\\
&&+ \frac{3B_{G}^4L_{F}^2 \Delta}{2a\alpha} \mathbb{E} \| x_{t - \tau_{t}^{s}}^{s} - x_{t}^{s}\|^2
\nonumber\\
&\overset{(\ref{young's})}{\leq}& (\frac{\alpha}{2} + \frac{3B_{G}^4L_{F}^2 }{2a\alpha})\mathbb{E} \|x_{t}^{s} - x^{*}\|^2 + \frac{3B_{G}^4L_{F}^2 }{2a\alpha} \mathbb{E} \|\widetilde{x}^s  - x^{*}\|^2 
\nonumber\\
&&+ \frac{3 \eta^2 B_{G}^4L_{F}^2 \Delta T}{2a\alpha}  \sum_{l=t-\tau_{t}^{s}}^{t-1} \mathbb{E} \| \nabla \widehat{f}(x_{l}^{s}) \|^2.
\end{eqnarray}
We then bound $T_6$ with
\begin{eqnarray} \label{distribute_T9}
-T_6 &=& \mathbb{E} \langle x_{t}^{s} - x^{*}, \nabla f_{i_t j_{t}}(x_{t - \tau_{t}^{s}}^{s}) \rangle
\nonumber\\
&=& \underbrace{ \mathbb{E} \langle x_{t - \tau_{t}^{s}}^{s} - x^{*}, \nabla f_{i_t j_{t}}(x_{t - \tau_{t}^{s}}^{s}) \rangle }_{T_{8}}+ \underbrace{ \sum_{l=t - \tau_{t}^{s}}^{t-1} \mathbb{E} \langle x_{l+1}^{s} - x_{l}^{s}, \nabla f_{i_t j_{t}}(x_{l}^{s}) \rangle }_{T_{9}}
\nonumber\\
& & + \underbrace{\sum_{l=t - \tau_{t}^{s}}^{t-1} \mathbb{E} \langle x_{l+1}^{s} - x_{l}^{s},\nabla f_{i_t j_t}(x_{t - \tau_{t}^{s}}^{s}) - \nabla f_{i_t j_t}(x_l^{s}) \rangle }_{T_{10}}.
\end{eqnarray}
Because $f_{ij}(x)$ is convex, we can bound $T_{8}$ by
\begin{equation} \label{distribute_T11}
T_{8} \geq \mathbb{E}[f_{i_t j_{t}}(x_{t - \tau_{t}^{s}}^{s}) - f_{i_t j_{t}}(x^{*})] = \mathbb{E}f(x_{t - \tau_{t}^{s}}^{s}) - f(x^{*}).
\end{equation}
We then bound $T_{9}$ by
\begin{eqnarray} \label{distribute_T12}
T_{9} &\geq& \sum_{l=t - \tau_{t}^{s}}^{t-1} \mathbb{E} [f_{i_t j_t}(x_{l+1}^{s}) - f_{i_t j_t}(x_{l}^{s}) - \frac{L_{f}}{2}\|x_{l}^{s} - x_{l+1}^{s}\|_{i_t j_t}^2]
\nonumber\\
& \geq& \mathbb{E}[ f_{i_t j_t}(x_{t}^{s}) - f_{i_t j_t}(x_{t - \tau_{t}^{s}}^{s}) - \frac{L_{f} \Delta}{2} \sum_{l=t - \tau_{t}^{s}}^{t-1} \|x_{l}^{s} - x_{l+1}^{s}\|^2]
\nonumber\\
&=& \mathbb{E}(f(x_{t}^{s}) - f(x_{t - \tau_{t}^{s}}^{s})) - \frac{L_{f} \Delta}{2} \sum_{l=t - \tau_{t}^{s}}^{t-1} \mathbb{E} \|x_{l}^{s} - x_{l+1}^{s}\|^2,
\end{eqnarray}
where $\| \cdot\|_{i_t j_t}$ denotes the support of $\nabla f_{i_t j_t}(x)$, and the second inequality comes from Assumption \ref{assumption:sparsity}.
We proceed to bound $T_{10}$ by
\begin{eqnarray} \label{distribute_T13}
-T_{10} &=& \sum_{l=t - \tau_{t}^{s}}^{t-1} \mathbb{E} \langle x_{l}^{s} - x_{l+1}^{s}, \nabla f_{i_t j_t}(x_{t - \tau_{t}^{s}}^{s}) - \nabla f_{i_t j_t}(x_{l}^{s}) \rangle
\nonumber\\
& \overset{(\ref{cauthy})}{\leq} & \sum_{l=t - \tau_{t}^{s}}^{t-1} \mathbb{E} \|x_{l}^{s} - x_{l+1}^{s}\|_{i_t j_t} \|\nabla f_{i_t j_t}(x_{t - \tau_{t}^{s}}^{s}) - \nabla f_{i_t j_t}(x_{l}^{s})\|_{i_t j_t}
\nonumber\\
& \leq & \sum_{l=t - \tau_{t}^{s}}^{t-1}  \mathbb{E} \|x_{l}^{s} - x_{l+1}^{s}\|_{i_t j_t} \sum_{j=t - \tau_{t}^{s}}^{l-1} \mathbb{E} \|\nabla f_{i_t j_t}(x_{j}^{s}) - \nabla f_{i_t j_t}(x_{j+1}^{s})\|_{i_t j_t}
\nonumber\\
&\leq& L_{f} \sum_{l=t - \tau_{t}^{s}}^{t-1} \mathbb{E} \|x_{l}^{s} - x_{l+1}^{s}\|_{i_t j_t} \sum_{j=t - \tau_{t}^{s}}^{l-1} \mathbb{E} \|x_{j}^{s} - x_{j+1}^{s}\|_{i_t j_t}
\nonumber\\
&\overset{(\ref{cauthy})}{\leq}& \frac{L_{f}}{2} \sum_{l=t - \tau_{t}^{s}}^{t-1} \sum_{j=t - \tau_{t}^{s}}^{l-1} \mathbb{E} ( \|x_{l}^{s} - x_{l+1}^{s}\|_{i_t j_t}^2 + \|x_{j}^{s} - x_{j+1}^{s}\|_{i_t j_t}^2)
\nonumber\\
&\leq&\frac{L_{f}}{2} (t-(t - \tau_{t}^{s})-1) \sum_{l=t - \tau_{t}^{s}}^{t-1} \mathbb{E}\|x_{l}^{s} - x_{l+1}^{s}\|_{i_t j_t}^2
\nonumber\\
&\leq& \frac{L_{f}}{2} (T -1) \sum_{l=t - \tau_{t}^{s}}^{t-1} \mathbb{E} \|x_{l}^{s} - x_{l+1}^{s}\|_{i_t j_t}^2
\nonumber\\
&\leq& \frac{\Delta L_{f} (T - 1) \eta^2}{2} \sum_{l=t - \tau_{t}^{s}}^{t-1} \mathbb{E} \| \nabla \widehat{f}(x_{l}^{s}) \|^2,
\end{eqnarray}
where the second inequality comes from the triangle inequality, and the fifth inequality can be obtained by using a simple counting argument, and the last inequality comes from Assumption \ref{assumption:sparsity}.
Substituting (\ref{distribute_T11}), (\ref{distribute_T12}), (\ref{distribute_T13}) into (\ref{distribute_T9}), we have
\begin{equation}\label{distribute_T9_2}
T_6 \leq \mathbb{E}(f(x^{*}) - f(x_{t}^{s})) + \frac{\Delta L_{f} T \eta^2}{2} \sum_{l=t - \tau_{t}^{s}}^{t-1} \mathbb{E} \| \nabla \widehat{f}(x_{l}^{s}) \|^2.
\end{equation}
Then we substitute (\ref{distribute_T8_2}) and (\ref{distribute_T9_2}) into (\ref{distribute_T7}) so as to obtain
\begin{eqnarray} \label{distribute_T7_2}
T_4 &\leq& f(x^{*}) - f(x_{t}^{s}) + (\frac{\alpha}{2} + \frac{3  B_{G}^4 L_{F}^2}{2a\alpha}) \mathbb{E} \|x_{t}^{s} - x^{*}\| + \frac{3 B_{G}^4 L_{F}^2}{2a\alpha} \mathbb{E} \|\widetilde{x}^s - x^{*}\|^2 
\nonumber\\
& &+ (\frac{\Delta L_f T \eta^2}{2} + \frac{3 B_{G}^4 L_{F}^2 \Delta T \eta^2}{2a\alpha}) \sum_{l=t - \tau_{t}^{s}}^{t-1} \mathbb{E} \| \nabla \widehat{f}(x_{l}^{s}) \|^2.
\end{eqnarray}
Substituting (\ref{distribute_T7_2}) into (\ref{Theorem_2}), we have
\begin{eqnarray}\label{Theorem_2_2}
&&\mathbb{E} \|x_{t+1}^{s} - x^{*}\|^2 
\nonumber\\
&\leq& \mathbb{E} \|x_{t}^{s} - x^{*}\|^2 + \eta^2 \mathbb{E} \|\nabla \widehat{f}(x_{t - \tau_{t}^{s}}^{s})\|^2 + 2\eta (\frac{\alpha}{2} + \frac{3B_{G}^4 L_{F}^2 }{2a\alpha}) \mathbb{E}\|x_{t}^{s} - x^{*}\|^2
\nonumber\\
& & + \frac{3\eta B_{G}^4 L_{F}^2 \Delta}{a\alpha} \mathbb{E} \|\widetilde{x}^s - x^{*}\|^2 + (\Delta L_f T \eta^3 + \frac{3 B_{G}^4 L_{F}^2 \Delta T \eta^3}{a\alpha}) \sum_{l=t - \tau_{t}^{s}}^{t-1} \mathbb{E} \|\nabla \widehat{f}(x_{l}^{s} )\|^2
\nonumber\\
& & - 2\eta \mathbb{E} (f(x_{t}^{s}) - f(x^{*})).
\end{eqnarray}
Summing up this inequality from $t=0$ to $t=K-1$, we get
\begin{eqnarray}\label{Theorem_1_3}
&&\mathbb{E} \|x_{K}^{s} - x^{*}\|^2 
\nonumber\\
&\leq& \mathbb{E} \|\widetilde{x}^{s} - x^{*}\|^2 + \eta^2 \sum_{t=0}^{K-1} \mathbb{E} \|\nabla \widehat{f}(x_{t - \tau_{t}^{s}}^{s})\|^2 + 2\eta (\frac{\alpha}{2} + \frac{3B_{G}^4 L_{F}^2 }{2a\alpha})\sum_{t=0}^{K-1} \mathbb{E} \|x_{t}^{s} - x^{*}\|^2
\nonumber\\
& & + \frac{3\eta B_{G}^4 L_{F}^2 }{a\alpha} K \mathbb{E} \|\widetilde{x}^s - x^{*}\|^2 + (\Delta L_f T \eta^3 + \frac{3 B_{G}^4 L_{F}^2 \Delta \tau \eta^3}{a\alpha}) \sum_{t=0}^{K-1} \sum_{l=t-\tau_{t}^{s}}^{t-1} \mathbb{E} \|\nabla \widehat{f}(x_{l}^{s})\|^2
\nonumber\\
& & - 2\eta \sum_{t=0}^{K-1} \mathbb{E} (f(x_{t}^{s}) - f(x^{*}))
\nonumber\\
&\leq&(1 + \frac{3 \eta B_{G}^4 L_{F}^2 K}{a\alpha}) \mathbb{E} \|\widetilde{x}^s - x^{*}\|^2 + (\eta \alpha + \frac{3\eta B_{G}^4 L_{F}^2 }{a \alpha}) \sum_{t=0}^{K-1} \mathbb{E} \|x_{t}^{s} - x^{*}\|^2
\nonumber\\
& & +[\eta^2 + \Delta L_{f}T^2 \eta^3 + \frac{3 \eta^3 B_{G}^4 L_{F}^2 \Delta T^2}{a\alpha}] \sum_{t=0}^{K-1} \mathbb{E} \|\nabla \widehat{f}(x_{t - \tau_{t}^{s}}^{s})\|^2
\nonumber\\
&& - 2\eta \sum_{t=0}^{K-1} \mathbb{E} (f(x_{t}^{s}) - f(x^{*}))
\nonumber\\
&\leq& (1 + \frac{3 \eta B_{G}^4 L_{F}^2 K}{a\alpha}) \frac{2}{\mu_{f}} \mathbb{E} (f(\widetilde{x}^s) - f(x^{*}) ) 
\nonumber\\
&&- (2\eta - (\eta \alpha + \frac{3\eta B_{G}^4 L_{F}^2 }{a \alpha}) \frac{2}{\mu_f} ) \sum_{t=0}^{K-1} \mathbb{E} ( f(x_{t}^{s}) - f(x^{*}) )
\nonumber\\
& & +[\eta^2 + \Delta L_{f}T^2 \eta^3 + \frac{3 \eta^3 B_{G}^4 L_{F}^2 \Delta T^2}{a\alpha}]  \sum_{t=0}^{K-1} \mathbb{E} \|\nabla \widehat{f}(x_{t - \tau_{t}^{s}}^{s})\|^2.
\end{eqnarray}
Combining Lemma~\ref{lemma:5} and Lemma~\ref{lemma:1} and setting $\alpha = \frac{\mu_{f}}{8}$, we can get
\begin{equation}\label{distribute_last}
\frac{1}{K} \sum_{t=0}^{K-1} \mathbb{E} ( f(x_{t}^{s}) - f(x^{*}) ) \leq \frac{\frac{\mu_{f}}{2}  + PQRK + U}{\frac{7}{4} \eta K - PQRK - U} \mathbb{E}(f(\widetilde{x}^s) - f(x^{*})).
\end{equation}
where $U, P, Q$ and $R$ are defined in (\ref{theorem_distribution_1:parameter}).
At last, taking (\ref{jeason}) into (\ref{distribute_last}), we can complete the proof by
\begin{equation}\label{distribute_last_2}
\mathbb{E} ( f(\widetilde{x}^{s+1}) - f(x^{*}) ) \leq \frac{\frac{\mu_{f}}{2}  + PQRK + U}{\frac{7}{4} \eta K - PQRK - U} \mathbb{E}(f(\widetilde{x}^s) - f(x^{*})).
\end{equation}

\noindent \textbf{Proof of Corollary 2.}

\noindent The proof of Corollary 2 is analogous to that of Corollary 1.

\bibliographystyle{spbasic}      
\bibliography{Reference}   

\end{document}